\documentclass[10pt,a4paper]{article}

\usepackage{amsthm}
\usepackage{amsmath}
\usepackage{amssymb}
\usepackage{empheq}
\usepackage{stackrel}

\newcommand{\ds}{\displaystyle}

\def\arccosh{\operatorname{arccosh}}
\def\Vol{\operatorname{Vol}}
\def \CC{\mathbb C}
\def \Re{\operatorname{Re}}
\def \arg{\operatorname{arg}}
\newtheorem{conjecture}{Conjecture}
\newtheorem{theorem}{Theorem}
\newtheorem{lemma}{Lemma}
\newtheorem{proposition}{Proposition}
\newtheorem{remark}{Remark}

\begin{document}

\title{Asymptotic Behavior of Colored HOMFLY Polynomial of Figure Eight Knot}

\author{Ka Ho Wong, Thomas Kwok-Keung Au}

\date{}

\maketitle

\begin{abstract}
In this paper we investigate the asymptotic behavior of the colored HOMFLY polynomial of the figure eight knot associated with the symmetric representation.  We establish an analogous asymptotic expansion for the colored HOMFLY polynomial.   From the asymptotic behavior we show that the Chern-Simons invariants and twisted Reidemeister torsion can be obtained with suitable modification of the case of colored Jones polynomial.
\end{abstract}

\section{Introduction}
This paper aims to find out what kinds of information can be extracted from the asymptotic behavior of the colored HOMFLY polynomial for a knot.   Our study starts from an understanding of the known asymptotic behavior of the colored Jones polynomial.  Its related historical background is briefly described below; mainly summarized in the works of H.Murakami.  From the known results and the development, the study of this paper is naturally motivated and a description is given below.  Our main theorem is then stated and the study is outlined.

The asymptotic behavior of the colored Jones polynomial, or $SU(2)$ invariant, has been investigated for a very long time. It started from the classical volume conjecture (Conjecture~\ref{cvc} below), which says that the evaluation of colored Jones polynomial of a knot $K$ at an $N$-th root of unity captures the simplicial volume of the knot complement $\mathbb{S}^{3} \backslash K$.

\begin{conjecture}\label{cvc}(Classical volume conjecture \cite{K97,MM01}) Let $K$ be a knot and $J_{N}^{(2)}(K;q)$ be the colored Jones polynomial of $K$ evaluated at $q$. We have
\begin{align*}
\lim_{N \to \infty} \dfrac{\log|J_{N}(K;e^{\frac{2\pi i}{N}})|}{N} = \frac{\operatorname{Vol}(\mathbb{S}^{3} \backslash K)}{2\pi},
\end{align*}
where $\operatorname{Vol}(\mathbb{S}^{3} \backslash K)$ is the simplicial volume of the knot complement.
\end{conjecture}

Several generalizations of the volume conjecture have been proposed, for example, see \cite{HM11} for a general review and \cite{DG11} for the physical interpretation.  A particular example is that the asymptotic behavior of the colored Jones polynomial captures the Chern-Simons invariant together with the Reidemeister torsion of the knot. A special case of the conjecture has been proved by H.Murakami in \cite{HM13}.

\begin{theorem}\label{asymsu2}
{\em (Asymptotic expansion for $SU(2)$ invariant of $4_{1} $\cite{HM13})\/}
Let $u$ be a real number with $0 < u < \log((3+\sqrt{5})/2) = 0.9624\dots$ and put $\xi = 2\pi i + u$. Then we have the following asymptotic equivalence of the $SU(2)$ invariant of the figure-eight knot~$4_{1}$:
\begin{align*}
J^{(2)}_{N}(4_{1} ; \exp(\xi / (N+n-2)) \stackrel[N \to \infty]{\sim}{ }  \frac{\sqrt{\pi}}{2 \sinh (u/2)} T(u)^{1/2} \left( \frac{N}{\xi}\right)^{1/2} \exp \left(\frac{N}{\xi} S(u)\right),
\end{align*}
where
$$ S(u) = \operatorname{Li}_{2}\left(e^{u-\varphi(u)}\right) - \operatorname{Li}_{2}\left(e^{u+\varphi(u)}\right) - u\varphi(u) $$
and
$$ T(u) = \frac{2}{\sqrt{(e^{u} + e^{-u} +1)(e^{u} + e^{-u} -3)}}.$$
Here $\ds \varphi(u) = \arccosh(\cosh(u) - 1/2)$ and
$$ \operatorname{Li}_{2}(z) = - \int_{0}^{z} \frac{ \log(1-x)}{x} dx$$
is the dilogarithm function.
\end{theorem}

\subsection*{Motivation}

Although the asymptotic behavior of colored Jones polynomial draws a lot of attention to mathematicians, the asymptotic behavior of its generalization, colored HOMFLY polynomial, or $SU(n)$ invariant, does not. One reason is that the explicit formula for the colored HOMFLY polynomial is only known for a few knots. Fortunately for the figure eight knot we know much more. In particular, in~\cite{CLZ16} the classical volume conjecture has been extended to the colored HOMFLY polynomial associated with the symmetric representation, which is as follows.

\begin{conjecture}(Volume conjecture for $SU(n)$ invariant)\label{cvcn}
Let $K$ be a hyperbolic knot and $J_{N}^{(n)}(K;q)$ be the colored HOMFLY polynomial, or $SU(n)$ invariant, of $K$ associated with the symmetric representation evaluated at $q$. For $a=0,1, 2, \dots, n-2, s \in \mathbb{Z}$, we have
$$ 2\pi s\, \lim_{N \to \infty} \frac{\log J_{N}^{(n)}\left(K; \exp\left(\frac{2s\pi i}{N+a}\right)\right)}{N} = \operatorname{Vol}(\mathbb{S}^{3}\backslash K) + i\operatorname{CS}(\mathbb{S}^{3}\backslash K)$$
\end{conjecture}

In \cite{CLZ16} the conjecture is proved for the figure eight knot~$4_1$. Therefore, it is natural to aim for generalizing Theorem~\ref{asymsu2} to the colored HOMFLY polynomial. One crucial question is that what kind of information can be obtained from the asymptotic formula. Since the colored Jones polynomial, or $SU(2)$ invariant, captures the Reidemeister torsion associated with the $SL(2;\mathbb{C})$ representation of the knot group, it is natural to guess that the colored HOMFLY polynomial, or $SU(n)$ invariant, should capture the higher dimensional Reidemeister torsion. The higher dimensional Reidemeister torsion of knot complement has been explored by several authors (see for example \cite{PJ14}) and the torsion itself has a very interesting property relating to the volume of the knot complement.

\begin{theorem}{\em(\cite{PJ14})\/}
Let M be a connected, complete, hyperbolic 3-manifold of finite volume. Denote $\tau_{n}$ to be the $n$ dimensional Reidemeister torsion of M. Then
\begin{align*}
\lim_{k \to \infty} \dfrac{\log | \tau_{2k+1} (M)|}{(2k+1)^{2}} = - \dfrac{\Vol(M)}{4\pi}
\end{align*}
In addition, if $\eta$ is an acyclic spin structure on $M$, then
\begin{align*}
\lim_{k \to \infty} \dfrac{\log | \tau_{2k} (M;\eta)|}{(2k)^{2}} = - \dfrac{\Vol(M)}{4\pi}
\end{align*}
\end{theorem}

Combining the above observations, it is exciting to see whether the above theorem could be placed into the context of asymptotic expansion of $SU(n)$ invariant.

\subsection*{Main Result}

To test the validity of the idea given above, the first thing is to find out explicitly the asymptotic behavior of $SU(n)$ invariant.  Following similar ideas as in \cite{HM13}, we obtain the main result of this paper stated below.

\begin{theorem}\label{mainthm}
{\em (Asymptotic expansion for $SU(n)$ invariant of $4_{1}$)\/}
For even $\ds n \geq 2$, let u be a real number with $0 < u < \log((3+\sqrt{5})/2) = 0.9624\dots$ and put $\xi = 2\pi i + u$. Then we have the following asymptotic equivalence of the SU(n) invariant of the figure-eight knot $4_{1}$:
\begin{align}\label{asymsun}
&\lefteqn{ J^{(n)}_{N}(4_{1} ; \exp(\xi / (N+n-2)) \quad\stackrel[N \to \infty]{\sim}{ }  ((1-e^{u-\phi(u)})(1-e^{\phi(u)}))^{n-2} \times} \notag \\
&\frac{1}{(n-2)!}\frac{1}{(e^{u}-1)^{n-2}}\frac{\sqrt{-\pi}}{2 \sinh (u/2)}T(u)^{1/2}(\frac{N+n-2}{\xi})^{\frac{1}{2} +(n-2)} \exp (\frac{N+n-2}{\xi} S(u)),
\end{align}
where $S(u)$, $T(u)$ and $\phi(u)$ are defined as in Theorem~\ref{asymsu2}.
\end{theorem}

\subsection*{Plan of this paper}

The first part of this paper is to show the main theorem. In Section~\ref{sec2}, we will outline the proof the main theorem.  In the process, a number of propositions and lemmas will be stated only for clarity of the strategy.  The detailed proofs of these propositions and lemmas are delayed and collected in Section~\ref{sec3}.

In Section~\ref{sec4}, we discuss why the same method cannot apply to other roots of unity. Finally in Section~\ref{sec5}, we discuss some difficulties we meet in the process of verifying the conjecture we mentioned in the introduction.  This may highlight possible paths for further developing the theory.

\section{Proof Outline of the Main Theorem}\label{sec2}

In this section we borrow the idea in~\cite{HM13} to find out the asymptotic expansion formula. First of all, the $SU(n)$ invariant of figure eight knot~$4_1$ associated with symmetric representation is given by (\cite{CLZ16}, \cite{IMM12}), namely,
$$ J_{N}^{(n)}(4_{1};q) = \frac{1}{[n-2]!} \sum_{k=0}^{N-1} \frac{[n-2+k]!}{[k]!} q^{-k\left(N+\frac{n-2}{2}\right)} \prod_{l=1}^{k}\left(1-q^{N-l}\right)\left(1-q^{N+l+n-2}\right)\,. $$

Here we use the convention that $[n] = q^{n/2} - q^{-n/2}$. This $ J_{N}^{(n)}$ is reduced to the colored Jones polynomial by putting $n=2$.  Secondly,
recall the definition of quantum dilogarithm $S_{\gamma}(z)$ (see \cite{HM13,F95}); that is,
$$ S_{\gamma}(z)= \exp\left(\frac{1}{4}\int_{C_{R}}\frac{e^{zt}}{\sinh(\pi t)\sinh(\gamma t)}\frac{dt}{t}\right) \,, $$
where $|\operatorname{Re}(z)|< \pi + \operatorname{Re}(\gamma)$ and $C_{R}$ is the contour $(-\infty,-R]\cup \Omega_{R} \cup [R, \infty)$ with a semi-circle $\Omega_{R} = \{ R\exp (i(\pi - s)) | 0\leq s \leq \pi \}$ for $0<R<\min\{\pi / |\gamma|,1\}$. The poles of the integrand are $0,\pm i, \pm 2i, \dots$ and $\pm \pi i/\gamma, \pm 2 \pi i \gamma, \dots$.

The following formula from \cite[Lemma~2.2]{HM13} is very helpful to rewrite the quantum dilogarithm~$S_\gamma$ as exponents or vice versa.  It will be used frequently in our calculations.  The proof can be referred in the cited paper.

\begin{lemma}\label{Slemma} If $| \operatorname{Re}(z) |  < \pi $, then we have
\begin{align}
(1+e^{iz})S_{\gamma}(z+\gamma)=S_{\gamma}(z-\gamma)
\end{align}
\end{lemma}
Using this formula, we may rewrite the $SU(n)$ invariant of the figure eight knot, $J_N^{(n)}\left(4_1; q\right)$ in terms of quantum dilogarithms for $q=\exp\left(\frac{\xi}{N+n-2}\right)$.  First,
applying Lemma~\ref{Slemma} with the values
$$\ds \gamma=\frac{2\pi - iu}{2(N+n-2)}\,,\quad\xi = 2\pi i + u \quad\text{and}\quad z=\pi - iu -2(n-2+l)\gamma $$
and observing that $\ds \frac{\xi}{N+n-2}=2i\gamma $, we have
\begin{align}\label{F1}
\prod_{l=1}^{k}\left(1-e^{\frac{N-l}{N+n-2}\xi}\right) = \frac{S_{\gamma}(\pi - iu - (2(n+k-2)+1)\gamma)}{S_{\gamma}(\pi - iu - (2n-3)\gamma)}
\end{align}
Similarly, putting $z=-\pi -iu +2l\gamma$, we have
\begin{align}\label{F2}
\prod_{l=1}^{k}\left(1-e^{\frac{N+l+n-2}{N+n-2}\xi}\right) = \frac{S_{\gamma}(-\pi - iu + \gamma)}{S_{\gamma}(-\pi - iu + (2k+1)\gamma)}
\end{align}
On the other hand,
\begin{align*}
\frac{[n-2+k]!}{[k]!}
&= [n-2+k][n-2+k-1]\dots [k+2][k+1] \\
&=q^{-\frac{n-2}{2}\left(\frac{n+2k-1}{2}\right)}\prod_{l=k+1}^{n-2+k}(1-q^{l})
\end{align*}
Putting $\ds q=\exp \left(\frac{\xi}{N+n-2}\right)$, we have
\begin{align*}
1-q^{l}=1-e^{\frac{\xi l}{N+n-2}}
= 1 + e^{i(-\pi+2l\gamma)}
\end{align*}
By Lemma~\ref{Slemma} with $z=-\pi  +2l\gamma$, the product terms cancel each other successively,
\begin{align}
\prod_{l=k+1}^{n-2+k}\left(1 + e^{i(-\pi+2l\gamma)}\right)
&=\prod_{l=k+1}^{n-2+k} \frac{S_{\gamma}(-\pi  +(2l-1)\gamma)}{S_{\gamma}(-\pi + (2l+1)\gamma)} \notag\\ &=\frac{S_{\gamma}(-\pi + (2(k+1)-1)\gamma)}{S_{\gamma}(-\pi + (2(n+k-2)+1)\gamma)} \label{F3}
\end{align}
By (\ref{F1}), (\ref{F2}) and (\ref{F3}), the $SU(n)$~invariant is expressed in terms of quantum dilogarithms,
\begin{align*}
J_{N}^{(n)}\left(4_{1}, e^{\frac{\xi}{N+n-2}}\right)
&= \frac{1}{[n-2]!}\frac{S_{\gamma}(-\pi-iu+\gamma)}{S_{\gamma}(\pi-iu-(2n-3)\gamma)} \times \hbox{\hspace*{5em}} \\
\sum_{k=0}^{N-1}e^{-ku-\tfrac{(n-2)(n-1)\xi}{4(N+n-2)}}
&\frac{S_{\gamma}(\pi-iu-(2n+2k-3)\gamma)\,S_{\gamma}(-\pi+(2k+1)\gamma)}{S_{\gamma}(-\pi-iu+(2k+1)\gamma)\,S_{\gamma}(-\pi+(2n+2k-3)\gamma)}\,.
\end{align*}

In order to obtain the asymptotic expansion for the above, we need to rewrite the summation terms into a contour integral so that an estimate can be achieved.  For that purpose, define
\begin{align*}
g_{N+n-2}(z) = e^{-(N+n-2)uz}\frac{S_{\gamma}(\pi - iu + i(z + \frac{n-2}{N+n-2})\xi)S_{\gamma}(-\pi - iz\xi)}{S_{\gamma}(-\pi - iu - iz\xi)S_{\gamma}(-\pi - i(z+\frac{n-2}{N+n-2})\xi)}
\end{align*}
Since $S_{\gamma}(z)$ is defined for $|\operatorname{Re}(z)|<\pi + \operatorname{Re}(\gamma)$ and $\operatorname{Re}(\gamma)>0$, one may check that $g$ is well-defined if $z \in D$ where
\begin{align*}
D &=\left\{
x+iy \in \mathbb{C}~\left|
\begin{array}{c}
-\frac{2\pi x}{u} -\frac{\operatorname{Re}(\gamma)}{u}< y < \frac{2\pi}{u}(1-x)+\frac{\operatorname{Re}(\gamma)}{u} \qquad\text{and} \\
- \frac{2\pi}{u} (x+\frac{n-2}{N+n-2}) -\frac{\operatorname{Re}(\gamma)}{u}< y < \frac{2\pi}{u}(1-(x+\frac{n-2}{N+n-2})) +\frac{\operatorname{Re}(\gamma)}{u} \\
\end{array}\right.\right\}\\
&\\
&=\left\{
  x+iy \in \mathbb{C}~\left| -\frac{2\pi x}{u} -\frac{\operatorname{Re}(\gamma)}{u}< y <\frac{2\pi}{u}(1-(x+\frac{n-2}{N+n-2})) +\frac{\operatorname{Re}(\gamma)}{u} \right.\right\}
\end{align*}
Next, for small $\ds \epsilon > \frac{n-2}{N+n-2}$, define the contour $C(\epsilon) = C_{+}(\epsilon)\cup C_{-}(\epsilon)$ with the polygonal lines $C_{\pm}(\epsilon)$ defined  by
\begin{align*}
C_{+}(\epsilon)&:\hspace{3em} 1-\epsilon \rightarrow 1- \frac{u}{2\pi} -\epsilon +i \rightarrow -\frac{u}{2\pi} + \epsilon +i \rightarrow \epsilon \\
C_{-}(\epsilon)&:\hspace{3em} \epsilon \rightarrow \epsilon + \frac{u}{2\pi} -i \rightarrow 1 - \epsilon+\frac{u}{2\pi} -i \rightarrow 1-\epsilon
\end{align*}
Note that all the singularity points of the function $z \mapsto\tan((N+n-2)\pi z$ are $\dfrac{2k+1}{2(N+n-2)} \in D$ for $k=0,1,2, \dots, N-1$.  Then, using Residue Theorem, we may express the $SU(n)$~invariant as
\begin{align*}
J_{N}^{(n)}(4_{1}, e^{\frac{\xi}{N+n-2}}) =& \frac{1}{[n-2]!}\frac{S_{\gamma}(-\pi-iu+\gamma)}{S_{\gamma}(\pi-iu-(2n-3)\gamma)} \times {} \\
& \frac{i e^{\frac{u}{2}} (N+n-2)}{2 e^{(\frac{n-2}{2})(\frac{n-1}{2})\frac{\xi}{N+n-2}}}
\int_{C(\epsilon)}\tan((N+n-2)\pi z)g_{N+n-2}(z)dz
\end{align*}
In order to estimate the integral, let
$$ G_{\pm}(N,n,\epsilon) = \displaystyle
\int_{C_{\pm}(\epsilon)}\tan((N+n-2)\pi z)g_{N+n-2}(z)dz\,.$$
As a result, one may rewrite
\begin{align*}
J_{N}^{(n)}\left(4_{1}, e^{\frac{\xi}{N+n-2}}\right)
&=  \frac{1}{[n-2]!}\frac{S_{\gamma}(-\pi-iu+\gamma)}{S_{\gamma}(\pi-iu-(2n-3)\gamma)} \,\times\, \\
&\qquad\qquad \frac{i e^{\frac{u}{2}} (N+n-2)}{2 e^{(\frac{n-2}{2})(\frac{n-1}{2})\frac{\xi}{N+n-2}}}(G_{+}(N,n, \epsilon) + G_{-}(N,n, \epsilon))
\end{align*}
The integral in $G_{\pm}$ may be splitted by adding and subtracting the same term as follows,
\begin{align*}
\lefteqn{ G_{\pm}(N,n,\epsilon) = }\qquad \\
&\pm i \int_{C_{\pm}(\epsilon)} g_{N+n-2}(z)dz + \int_{C_{\pm}(\epsilon)} (\tan((N+n-2)\pi z) \mp i)g_{N+n-2}(z) dz
\end{align*}
Intriguingly, according to the next proposition which will be proven in the upcoming section, the second integral term can be controlled and so decays asymptotically.
\begin{proposition}\label{tanapp}
There exists a constant $K_{1,\pm}$ independent of N and $\epsilon$ such that
$$ \left|   \int_{C_{\pm}(\epsilon)} (\tan((N+n-2)\pi z) \mp i)g_{N+n-2}(z) dz \right| < \frac{K_{1,\pm}}{N+n-2}\,.$$
\end{proposition}

Therefore, to arrive at the asymptotic expansion of $J_N^{(n)}$, it remains to approximate $g_{N+n-2}$.
First define a function
$$\Phi^{(n)}_{N}(z)=\frac{1}{\xi}\left[\operatorname{Li_{2}}\left(e^{u-\left(z+\frac{n-2}{N}\right)\xi}\right)+\operatorname{Li_{2}}\left(e^{z\xi}\right)
-\operatorname{Li_{2}}\left(e^{u+z\xi}\right) - \operatorname{Li_{2}}\left(e^{\left(z+\frac{n-2}{N}\right)\xi}\right) -uz\right]$$
Since $\operatorname{Li_{2}}$ is analytic in $\mathbb{C}\setminus [1,\infty)$, one may verify that the function $\Phi^{(n)}_{N}(z)$ is analytic in the region
$$ D' = \left\{ z= x+iy \in \mathbb{C} \left| -\frac{2\pi}{u}x < y < \frac{2\pi}{u}\left(1-\left(x+\frac{n-2}{N+n-2}\right)\right) \right.\right\} \subset D$$
Note that the contour $C(\epsilon)$ and the poles of $\tan((N+n-2)\pi z)$ lie inside $D'$.   The next statement provides a comparison between $g_{N+n-2}$ and $\Phi^{(n)}_{N+n-2}$.  Again, we delay the proof of the statement to the next section.
\begin{proposition}\label{exprepn}
Let $p(\epsilon)$ be any contour in the parallelogram bounded by $C(\epsilon)$ connecting from~$\epsilon$ to~$1-\epsilon$, then there exists a constant $K_{2}>0$ independent of $N$ and $\epsilon$ such that
\begin{align*}
\lefteqn{ \left| \int_{p(\epsilon)} g_{N+n-2}(z) dz - \int_{p(\epsilon)} \exp( (N+n-2) \Phi_{N+n-2}^{(n)}(z)) dz \right| \,\leq }\\
&\qquad\qquad \frac{K_{2} \log (N+n-2)}{N+n-2}
\max_{\omega \in p(\epsilon)} \left\{ \exp ((N+n-2) \operatorname{Re} \Phi_{N+n-2}^{(n)}(z) \right\}
\end{align*}
\end{proposition}
Since $\Phi_{N+n-2}^{(n)}(z)$ is analytic on $D'$, by Cauchy's theorem
\begin{align*}
&\int_{C_{+}(\epsilon)} \exp\left((N+n-2)\Phi_{N+n-2}^{(n)}(z)\right)dz \\
= &-\int_{C_{-}(\epsilon)} \exp\left((N+n-2)\Phi_{N+n-2}^{(n)}(z)\right)dz
\end{align*}
To approximate the above two integrals, we need the following generalized saddle point approximation, which will be proved in the next section.
\pagebreak
\begin{theorem}\label{FSA}{\em (One-parameter family version for saddle point approximation)}
Let $\{ \Phi_{y}(z)\}_{y \in [0,1]}$ be a family of holomorphic functions smoothly depending on $y \in [0,1]$. Let $C(y,t): [0,1]^{2} \to \CC$ be a continuous family of closed contours with length uniformly bounded above by a fixed constant $L$, such that for each $y \in [0,1]$, $C(y,t)$ lies inside the domain of $\Phi_{y}(z)$, for which $z_{y}$ is the only saddle point along the of contour $C_{y}$ and $\max\operatorname{Re}\left[\Phi_{y}(z)\right]$ is attained at $z_{y}$. Further assume that $\left|\arg\left( \sqrt{-\frac{d^{2}\Phi_{0}}{dz^{2}}(z_{0})}\right)\right| < \pi /4$. Then for each subsequence $\{y_{N} \}_{N \in \mathbb{N}}$ with $y_{N} \to 0$ as $N \to \infty$, we have the following generalized saddle point approximation:
\begin{align*}
\int_{C_{y_{N}}} \exp (N \Phi_{y_{N}}(z)) dz
=  \sqrt{\frac{2\pi}{N\left(-\frac{d^{2} \Phi_{y_{N}}}{dz^{2}} (z_{y_{N}})\right)}}\, \exp(N \Phi_{y_{N}}(z_{y_{N}})) \left(1 + O\left(\frac{1}{N}\right)\right)
\end{align*}
\end{theorem}
Applying Theorem~\ref{FSA} to our situation, we have
\begin{theorem}\label{saddlepointapp}{\em (Behavior of $\displaystyle\int_{C_{\pm}(\epsilon)} \exp\left(N\Phi_{N}^{(n)}(z)\right)dz$ for large $N$)\/}
Let $\ds z^{(n)}_{N}$ be the saddle point of $\ds \Phi^{(n)}_{N}$ inside the contour $C(\epsilon)$. Then
\begin{align*}
\int_{C_{-}(\epsilon)} \exp\left(N\Phi_{N}^{(n)}(z)\right)dz \quad\stackrel[N \to \infty]{\sim}{ }\quad  \frac{\sqrt{2\pi}\exp\left(N\Phi_{N}^{(n)}\left(z^{(n)}_{N}\right)\right)}{\sqrt{N}\sqrt{- \frac{d^{2} \Phi_{N}^{(n)}}{dz^{2}}\left(z^{(n)}_{N}\right)}}
\end{align*}
\end{theorem}
Together with the following proposition, which provides a control on the right-hand side, the integral in Theorem~\ref{saddlepointapp} is ensured to have exponentially growth.
\begin{proposition}\label{positive}
$\operatorname{Re}\Phi^{(n)}_{N+n-2}\left(z^{(n)}_{N}\right)$ is positive for $0<u<\log((3+\sqrt{5})/2)$.
\end{proposition}

Combining the controls in Propositions~\ref{tanapp} and~\ref{exprepn} and Theorem~\ref{saddlepointapp}, we are able to estimate $G_{\pm}(N,n,\epsilon)$, namely,
\begin{align*}
\lefteqn{\lim_{N \to \infty} \left| \frac{G_{\pm}(N+n-2,\epsilon)}{\displaystyle\pm i \int_{C_{\pm}(\epsilon)} \exp\left((N+n-2)\Phi^{(n)}_{N+n-2}(z)\right) dz} -1 \,\right| \quad\leq }\\
&\qquad\qquad \frac{K_{1, \pm}}{\displaystyle N\left| \int_{C_{\pm}(\epsilon)} \exp\left((N+n-2)\Phi^{(n)}_{N+n-2}(z)\right) dz \right|}\quad + \\
&\qquad\qquad\qquad
\frac{K_{2} \log(N+n-2)}{N+n-2} \times \frac{\exp\left((N+n-2) \operatorname{Re} \Phi^{(n)}_{N+n-2}(z^{(n)}_{N})\right)}{\displaystyle\left| \int_{C_{\pm}(\epsilon)} \exp\left((N+n-2)\Phi^{(n)}_{N+n-2})(z) \right) dz \right|} \\
&\qquad\qquad \xrightarrow{N \to \infty} 0\,.
\end{align*}
Thus, up to this point, we can asymptotically express $J_N^{(n)}$ in terms of quantum dilogarithm and a contour integral involving exponential of $N\Phi^{(n)}_N$.  That is,
\begin{align*}
\lefteqn{ J^{(n)}_{N}\left(4_1,e^{\xi/(N+n-2)}\right) \quad\stackrel[N \to \infty]{\sim}{ }\quad } \\
&\qquad\qquad \frac{1}{[n-2]!}\frac{S_{\gamma}(-\pi-iu+\gamma)}{S_{\gamma}(\pi-iu-(2n-3)\gamma)} \quad\times\quad \\
&\qquad\qquad\qquad\qquad \frac{(N+n-2)e^{u/2}}{2}\int_{C_{-}(\epsilon)}\exp\left((N+n-2)\Phi^{(n)}_{N+n-2}(z)\right)dz\,.
\end{align*}
Moreover, we also have the fact that (see \cite{HM13})
\begin{align*}
\lim_{N \to \infty} \frac{d^{2} \Phi_{N}^{(n)}}{dz^{2}}\left(z^{(n)}_{N}\right)
= \frac{d^{2} \Phi^{(2)}}{dz^{2}}\left(z^{(2)}\right)
= \xi \sqrt{(e^{u}+e^{-u}-1)(e^{u}+e^{-u}-3)}
\end{align*}

To obtain the whole asymptotic expansion of $J^{(n)}_N$ in term of $\Phi^{(n)}_N$, we need to study the asymptotic behavior of the quantum dilogarithm, as given in the lemma below.
\begin{lemma}\label{Sratio}
For $\gamma = \frac{2\pi - iu}{2(N+n-2)}$ with $u>0$ and an even integer $n$, we have
\begin{align*}
\frac{S_{\gamma}(-\pi-iu+\gamma)}{S_{\gamma}(\pi-iu-(2n-3)\gamma)}
&=
\frac{ e^{u\pi / \gamma}-1 }{ \prod_{k=0}^{n-2} \left(e^{u-2k\gamma i}-1\right) } \quad\stackrel[N \to \infty]{\sim}{ }\quad
\frac{e^{2\pi i u (N+n-2) / \xi}}{\left(e^{u}-1\right)^{n-1}}\,.
\end{align*}
\end{lemma}
In order to apply the saddle point approximation, we have to solve the equation
\begin{equation}\label{saddlepteqn}
\frac{d\Phi_{N+n-2}^{(n)}}{dz}(z)=0 \,.
\end{equation}
Recall that
\begin{align*}
\Phi^{(n)}_{N+n-2}(z) &= \\
&\hspace*{-2em} \frac{1}{\xi}\left[\operatorname{Li}_{2}\left(e^{u-\left(z+\frac{n-2}{N+n-2}\right)\xi}\right) + \operatorname{Li}_{2}\left(e^{z\xi}\right) - \operatorname{Li}_{2}\left(e^{u+z\xi}\right) - \operatorname{Li}_{2}\left(e^{\left(z+\frac{n-2}{N+n-2}\right)\xi}\right)\right]\, -uz \\
\frac{d}{d\mu}\operatorname{Li}_{2}(e^{\mu})
&= \operatorname{Li}_{1}(e^{\mu})=-\log(1-e^{\mu})
\end{align*}
The desired saddle point equation~(\ref{saddlepteqn}) can be rewritten as below,
\begin{align*}
\log\left(\frac{\left(1-e^{u-(z+\frac{n-2}{N+n-2})\xi}\right) \left(1-e^{u+z\xi}\right) \left(1-e^{(z+\frac{n-2}{N+n-2})\xi}\right)}{1-e^{z\xi}}\right)-u &=0 \,,\\
\intertext{which in turns becomes,}
 \frac{\left(1-e^{u-(z+\frac{n-2}{N+n-2})\xi}\right)\left(1-e^{u+z\xi}\right) \left(1-e^{(z+\frac{n-2}{N+n-2})\xi}\right)}{1-e^{z\xi}} &= e^{u}\,.
\end{align*}
With $\ds a=e^{u}$, $\ds b= e^{\frac{n-2}{N+n-2}\xi}$ and $\ds w = e^{z\xi}$, the above equation is equivalent to
\begin{align}\label{cubicpoly}
 ab^{2}w^{3} - (b^{2}+a^{2}b)w^{2} + (a^{2}+b)w -a =0
\end{align}
Let $w^{(n)}_{N+n-2}$ be the solution for~$w$ inside the domain $C(\epsilon)$ and $e^{z^{(n)}_{N+n-2}\xi}=w^{(n)}_{N+n-2}$. The asymptotic expansion formula of $J^{(n)}_N$ in Theorem~\ref{mainthm} is then obtained.

\begin{remark}
When $\ds n=2$ (i.e. $b=1$), after factoring out the factor $\ds (w-1)$ we obtained the quadratic equation appeared in \cite{HM13}. In this case $z^{(2)}_{N+n-2}= z^{(2)}$ is independent of $N$.
\end{remark}

The last step to establish Theorem~\ref{mainthm} is to change $\Phi^{(n)}_{N+n-2}$ into $\Phi^{(2)}$. The estimation between them is given by the following lemma, which is direct consequence of L'Hospital rule.

\begin{lemma}\label{diff1} For any $z \in D'$,
\begin{align*}
\lim_{N \to \infty} (N+n-2)\left(\Phi_{N+n-2}^{(n)}(z) - \Phi^{(2)}(z)\right) =  (n-2)\log\left((1-e^{u-z\xi})(1-e^{z\xi})\right)
\end{align*}
\end{lemma}

From Equation~(3.1) in \cite{HM13} we know that $z^{(2)}=\dfrac{\phi(u) + 2 \pi i}{ \xi}$. That means
\begin{align*}
\lefteqn{ \exp \left((N+n-2)\Phi^{(n)}_{N}\left(z^{(n)}_{N}\right)\right) } \\
&\stackrel[N \to \infty]{\sim}{}
\exp\left((n-2)\log((1-e^{u}(\omega^{(n)}_{N})^{-1})(1-\omega^{(n)}_{N})\right)\exp((N+n-2)(\Phi^{(2)}(z^{(n)}_{N}))) \\
&\stackrel[N \to \infty]{\sim}{} \exp\left((n-2)\log((1-e^{u}(\omega^{(2)})^{-1})(1-\omega^{(2)}))\right)\exp((N+n-2)(\Phi^{(2)}(z^{(n)}_{N})))\\
&\stackrel[N \to \infty]{\sim}{} \exp\left((n-2)\log((1-e^{u-\phi(u)})(1-e^{\phi(u)}))\right)\exp((N+n-2)(\Phi^{(2)}(z^{(n)}_{N})))\\
&\stackrel[N \to \infty]{\sim}{} \left((1-e^{u-\phi(u)})(1-e^{\phi(u)})\right)^{n-2}\exp((N+n-2)(\Phi^{(2)}(z^{(n)}_{N}))
\end{align*}

Using (\ref{cubicpoly}), one can show that $z^{(n)}_{N+n-2} - z^{(2)} = O\left(\frac{1}{N+n-2}\right)$. Together with the fact that $z^{(2)}$ satisfies the equation $\left.\dfrac{d\Phi^{(2)}}{d\omega}\right|_{z^{(2)}}=0$, we have
\begin{lemma}\label{diff2}
$\ds \lim_{N \to \infty} (N+n-2)\left(\Phi^{(2)}(z^{(n)}_{N+n-2}) - \Phi^{(2)}(z^{(2)})\right) = 0$
\end{lemma}

As a result,
\begin{align*}
\lefteqn{ \exp \left((N+n-2)\Phi^{(n)}_{N}(z^{(n)}_{N})\right) } \\
& \stackrel[N \to \infty]{\sim}{}\, \left((1-e^{u+\phi(u)})(1-e^{\phi(u)})\right)^{n-2}\exp((N+n-2)(\Phi^{(2)}(z^{(2)}))
\end{align*}

Finally we consider the large $N$ behavior of the term $\dfrac{1}{[n-2]!}$. Note that
\begin{align*}
\frac{1}{[k]}
= &\frac{1}{q^{-k/2}}\, \frac{1}{q^{k}-1}= \frac{1}{e^{-k\xi / 2(N+n-2)}}\, \frac{1}{e^{k\xi / (N+n-2)} -1} \\
= &\frac{1}{e^{-k\xi / 2(N+n-2)}}\, \frac{1}{(k\xi / (N+n-2))\left(\sum_{\ell=1}^{\infty}(k\xi / (N+n-2)^{\ell-1}/ \ell!\right)} \\
 \stackrel[N \to \infty]{\sim}{ } & \frac{N+n-2}{k\xi}
\end{align*}

Therefore by multiplying the terms together we get
\begin{align*}
\frac{1}{[n-2]!} \quad \stackrel[N \to \infty]{\sim}{ }\quad  \frac{1}{(n-2)!}\left(\frac{N+n-2}{\xi}\right)^{n-2}
\end{align*}

This complete the proof of Theorem~\ref{mainthm}.

\section{Proof of Results listed in Section~\ref{sec2}}\label{sec3}

\begin{proof}(Proof of Proposition~\ref{tanapp})
We follow the line of the proof in \cite{HM13} with suitable modification. Recall that for $| \operatorname{Re}(z) | < \pi$, or $| \operatorname{Re}(z)|= \pi$ and $\operatorname{Im}(z)>0$,
\begin{align*}
\frac{1}{2i}\operatorname{Li}_{2}(-e^{iz})
&=\frac{1}{4}\int_{C_{R}}\frac{e^{zt}}{t^{2}\sinh(\pi t)}dt \\
\text{\hskip -6em which implies \hskip 4em}
S_{\gamma}(z)
&= \exp\left(\frac{1}{2i\gamma}\operatorname{Li}_{2}(-e^{iz}) + I_{\gamma}(z)\right) \\
&= \exp\left(\frac{N+n-2}{\xi}\operatorname{Li}_{2}(-e^{iz}) + I_{\gamma}(z)\right), \\
\text{\hskip -7em where \hskip 7.5em}
I_{\gamma}(z) &= \frac{1}{4}\int_{C_{R}}\frac{e^{zt}}{t\sinh(\pi t)} \left(\frac{1}{\sinh(\gamma t)}-\frac{1}{\gamma t}\right)dt\,.
\end{align*}
Then the above $S_{\gamma}$ is substituted into the definition of $g_{N+n-2}$ and it leads to
\begin{align*}
\lefteqn{ g_{N+n-2}(z) = \exp\left[-(N+n-2)uz\right]\, \times }  \\
&\exp\left[\frac{N+n-2}{\xi}\left(\operatorname{Li}_{2}(e^{u-(z+\frac{n-2}{N+n-2})\xi}) + \operatorname{Li}_{2}(e^{z\xi}) - \operatorname{Li}_{2}(e^{u+z\xi})  - \operatorname{Li}_{2}(e^{(z+\frac{n-2}{N+n-2})\xi})\right)\right] \\
& \times \exp \left[ I_{\gamma}(\pi - iu + i(z + \frac{n-2}{N+n-2})\xi) + I_{\gamma}(-\pi - iz\xi)- I_{\gamma}(-\pi - iu - iz\xi) \right. \\
&\qquad  - \left. I_{\gamma}(-\pi - i\left(z+\frac{n-2}{N+n-2})\xi\right)\right]
\end{align*}
Let
\begin{align*}
\lefteqn{ \Phi^{(n)}_{N+n-2}(z) = } \qquad \\
&\frac{1}{\xi}(\operatorname{Li}_{2}(e^{u-(z+\frac{n-2}{N+n-2})\xi}) + \operatorname{Li}_{2}(e^{z\xi}) - \operatorname{Li}_{2}(e^{u+z\xi}) - \operatorname{Li}_{2}(e^{(z+\frac{n-2}{N+n-2})\xi}))-uz
\end{align*}
We have
\begin{align*}
\lefteqn{ g_{N+n-2}(z) = }\qquad \\
& \exp\left[(N+n-2)\Phi^{(n)}_{N+n-2}(z)\right]\, \times \\
& \hspace*{16pt} \exp \left[ I_{\gamma}\left(\pi - iu + i\left(z + \frac{n-2}{N+n-2}\right)\xi\right) + I_{\gamma}(-\pi - iz\xi) \right.\\
&\hspace*{36pt}  {} - \left. I_{\gamma}(-\pi - iu - iz\xi)  -I_{\gamma}\left(-\pi - i\left(z+\frac{n-2}{N+n-2}\right)\xi\right)\right]
\end{align*}

Decompose $C_{+}(\epsilon)$ as $C_{+,1}, C_{+,2}$ and $C_{+,3}$ by $\ds \epsilon \to (\epsilon -\frac{u}{2\pi} + i) \to (1 -\epsilon -\frac{u}{2\pi} + i ) \to 1-\epsilon$ and
$C_{-}(\epsilon)$ as $C_{-,1}, C_{-,2}$ and $C_{-,3}$ by $\ds \epsilon \to (\epsilon + \frac{u}{2\pi} - i) \to (1 -\epsilon +\frac{u}{2\pi} - i ) \to 1 - \epsilon$.\\
Write $I_{\pm,i}(N)$ be the integral along $C_{\pm,i}$ respectively. We are going to show the following controls on the integrals:
\begin{align}
|I_{+,1}(N+n-2)| &< \frac{K_{+,1}}{N+n-2} \label{prop2.11}\\
|I_{+,2}(N+n-2)| &< \frac{K_{+,2}}{N+n-2} \label{prop2.12}\ \\
|I_{+,3}(N+n-2)| &< \frac{K_{+,3}}{N+n-2} \label{prop2.13}\\\
|I_{-,1}(N+n-2)| &< \frac{K_{-,1}}{N+n-2} \label{prop2.14}\ \\
|I_{-,2}(N+n-2)| &< \frac{K_{-,2}}{N+n-2} \label{prop2.15}\ \\
|I_{-,3}(N+n-2)| &< \frac{K_{-,3}}{N+n-2} \label{prop2.16}\
\end{align}

Let us observe the comparison between $\Phi^{(n)}_{N}$ and $\Phi ^{(2)}$.  They are respectively related to the $SU(N)$ case and the $SU(2)$ case; with the latter one given in \cite{HM13}.
\begin{align*}
\Phi ^{(2)}(z) &= \frac{1}{\xi}(\operatorname{Li}_{2}(e^{u-\xi z}) - \operatorname{Li}_{2}(e^{u+\xi z})) - uz \\
\Phi^{(n)}_{N}(z)&=\frac{1}{\xi}\left(\operatorname{Li}_{2}(e^{u-(z+\frac{n-2}{N})\xi}) + \operatorname{Li}_{2}(e^{z\xi}) - \operatorname{Li}_{2}(e^{u+z\xi}) - \operatorname{Li}_{2}(e^{(z+\frac{n-2}{N})\xi})\right)-uz
\end{align*}
The proof of the above estimates for the contour integrals is basically the same as the one of Proposition~3.1 in \cite{HM13}.

 To prove (\ref{prop2.11}), first we estimate $|\tan((N+n-2)\pi((-u/2\pi + i)t + \epsilon)) - i|$. Note that
$$ |\tan((N+n-2)\pi((-u/2\pi + i)t + \epsilon)) - i|  \leq \frac{2e^{-2(N+n-2)\pi t}}{\left|e^{-2(N+n-2)\pi t - (N+n-2)uti + 2\epsilon i}\right|}  $$
Since $\epsilon$ can be arbitrary small as long as $N$ is large, for small $\epsilon$, by using (6.8) in \cite{HM13}, we have
$$ |\tan((N+n-2)\pi((-u/2\pi + i)t + \epsilon)) - i| \leq \frac{2e^{-2(N+n-2)\pi t}}{1-e^{-\pi^2 / u}} $$
So we have
$$ |I_{+,1}(N+n-2)| \leq \frac{2}{1-e^{-\pi^2 / u}} \int_{0}^{1}e^{-2N\pi t} \left|g_{N}((-\frac{u}{2\pi} +i)t + \epsilon)\right|$$
Recall the Lemma~6.1 in \cite{AH06} that for $|\operatorname{Re}(z)| \leq \pi$ we have
$$|I_{\gamma}(z)| \leq 2A + B|\gamma|\left(1+e^{-\operatorname{Im}(z)R}\right)$$
That means $\exp(\text{I part})$ is bounded above by some constant $M>0$ and
$$ \left|g_{N+n-2}((-\frac{u}{2\pi} +i)t + \epsilon)\right| \leq Me^{(N+n-2) \operatorname{Re}\Phi^{(n)}_{N}((\frac{u}{2\pi}+i)t+\epsilon)} $$
From the proof of~(6.2) in~\cite{HM13}, we know that $\operatorname{Re}\Phi^{(2)}((\frac{u}{2\pi}+i)t+\epsilon)<0$ for sufficiently small~$\epsilon>0$. Since $\Phi^{(n)}_{N+n-2} \longrightarrow \Phi^{(2)}$ as $N \to \infty$, we also have $\operatorname{Re}\Phi^{(n)}_{N+n-2}((\frac{u}{2\pi}+i)t+\epsilon)<0$ for $N$ large enough. Hence we have
$$\left|I_{+,1}(N+n-2)\right| \leq \frac{2}{1-e^{-\pi^{2} / u}}M \int_{0}^{1} e^{-2(N+n-2)\pi t} dt \leq \frac{K_{+,1}}{N+n-2}$$
This establishes the inequality~(\ref{prop2.11}).
The proof of the other inequalities~(\ref{prop2.12}--\ref{prop2.16}) are basically the same.
\end{proof}

\begin{proof}(Proof of Proposition~\ref{exprepn})
Write
\begin{align*}
g_{N}(z)
&= \exp((N+n-2)\Phi^{(n)}_{N+n-2}(z)) \times \exp (\text{I part})
\end{align*}

First, note that
\begin{align*}
&| \int_{p(\epsilon)} g_{N+n-2}(\omega) d\omega - \int_{p(\epsilon)} \exp( (N+n-2) \Phi_{N+n-2}^{(n)}(\omega)) d\omega | \\
= &| \int_{p(\epsilon)} \exp((N+n-2)\Phi_{N+n-2}^{(n)}) [ \exp(\text{I part}) - 1] | d\omega\\
\leq &\max_{\omega \in p(\epsilon)} \{ \exp ((N+n-2) \operatorname{Re} \Phi_{N+n-2}^{(n)}(\omega) \} \int_{p(\epsilon)} | \exp(\text{I part}) - 1| d\omega \\
= & \max_{\omega \in p(\epsilon)} \{ \exp ((N+n-2) \operatorname{Re} \Phi_{N+n-2}^{(n)}(\omega) \} \int_{\epsilon}^{1-\epsilon} |h_{\gamma}(\omega)| d\omega
\end{align*}

where
\begin{align*}
h_{\gamma}(\omega)= &\sum_{n=1}^{\infty} \frac{1}{n!}(I_{\gamma}(\pi - iu + i(\omega+\frac{n-2}{N+n-2})\xi)+I_{\gamma}(-\pi-i\omega\xi) \\
& -I_{\gamma}(-\pi - iu - i\omega\xi) - I_{\gamma}(-\pi-i(\omega+\frac{n-2}{N+n-2})\xi))^{n}
\end{align*}

In the above we use the analyticity of $h_{\gamma}(\omega)$ to change the contour to straight line parametrized by $t$, $t\in (\epsilon, 1-\epsilon)$.

Recall the lemma 3 in \cite{AH06} that there exist $A,B>0$ dependent only on $R$ such that if $|\operatorname{Re}(z)|<\pi$, we have
$$ |I_{\gamma}(z)| \leq A(\frac{1}{\pi-\operatorname{Re}(z)} + \frac{1}{\pi + \operatorname{Re}(z)})|\gamma| + B(1+e^{-\operatorname{Im}(z)R})|\gamma| $$

\pagebreak
So we have
\begin{enumerate}
\item
\begin{align*}
&|I_{1}|=|I_{\gamma}(\pi - iu + i\xi(t+\frac{n-2}{N+n-2}))| \\
\leq &A|\gamma| (\frac{1}{2\pi(t+\frac{n-2}{N+n-2})} + \frac{1}{2\pi - 2\pi(t+\frac{n-2}{N+n-2})}) \\
&+ B|\gamma|(1+e^{(u-u(t+\frac{n-2}{N+n-2}))R})\\
\leq &A|\gamma| (\frac{1}{2\pi(t+\frac{n-2}{N+n-2})} + \frac{1}{2\pi - 2\pi(t+\frac{n-2}{N+n-2})}) + B'|\gamma|
\end{align*}
\item
\begin{align*}
&|I_{2}|=|I_{\gamma}(-\pi - iu - i\xi t)|\\
\leq &A|\gamma| (\frac{1}{2\pi t} + \frac{1}{2\pi - 2\pi t}) + B|\gamma|(1+e^{(u-ut)R})\\
\leq &A|\gamma| (\frac{1}{2\pi t} + \frac{1}{2\pi - 2\pi t}) + B'|\gamma|
\end{align*}
\item
\begin{align*}
&|I_{3}|= |I_{\gamma}(-\pi  - i\xi t)|\\
\leq &A|\gamma| (\frac{1}{2\pi t} + \frac{1}{2\pi - 2\pi t}) + B|\gamma|(1+1)\\
\leq &A|\gamma| (\frac{1}{2\pi t} + \frac{1}{2\pi - 2\pi t}) + B'|\gamma|
\end{align*}
\item
\begin{align*}
&|I_{4}|= |I_{\gamma}(\pi - iu + i\xi(t+\frac{n-2}{N+n-2}))| \\
\leq &A|\gamma| (\frac{1}{2\pi(t+\frac{n-2}{N+n-2})} + \frac{1}{2\pi - 2\pi(t+\frac{n-2}{N+n-2})}) + B|\gamma|(1+1)\\
\leq &A|\gamma| (\frac{1}{2\pi(t+\frac{n-2}{N+n-2})} + \frac{1}{2\pi - 2\pi(t+\frac{n-2}{N+n-2})}) + B'|\gamma|
\end{align*}
\end{enumerate}

Let $f(t)=\dfrac{1}{t}+\dfrac{1}{1-t}$. Note that $f(t)\geq 4$ for $t \in [0,1]$.

From all four inequalities about $I_{\gamma}$, we have
\begin{align*}
|I_{1}+I_{2}-I_{3}-I_{4}| \leq |\gamma|(A'' f(t) + B' \frac{f(t)}{4}) \leq A''' |\gamma|f(t)
\end{align*}

Follow the argument in \cite{AH06}, p.537 we have
$\ds  \int_{|\gamma|}^{1-|\gamma|} f(t)^{n}dt \leq 2^{2n+1} \int_{|\gamma|}^{\frac{1}{2}} \frac{dt}{t^{n}} $ for $n \geq 1$.
Also since $|\gamma|=|\xi|/2(N+n-2)$ we have
\begin{align*}
\int_{|\gamma|}^{\frac{1}{2}} \frac{dt}{t} &= \log(N+n-2) - \log(|\gamma|) \leq \log(N+n-2) \\
\text{and }\int_{|\gamma|}^{\frac{1}{2}} \frac{dt}{t} &= \frac{1}{n-1}(\frac{1}{|\gamma|^{n-1}} - 2^{n-1}) \leq \frac{1}{|\gamma|^{n-1}} \text{ for $n \geq 2$}
\end{align*}

Therefore for $\epsilon > |\gamma|$ we have
\begin{align*}
\int_{\epsilon}^{1-\epsilon} |h_{\gamma}(t)| dt
\leq & \int_{|\gamma|}^{1-|\gamma|} |h_{\gamma}(t)| dt \\
\leq & \sum_{n=1}^{\infty} \frac{1}{n!} (A''')^{n} |\gamma|^{n} \int_{\epsilon}^{1-\epsilon} f(t)^{n} dt \\
\leq & 2|\gamma| (4A''' \log(N+n-2) + \sum_{n=2}^{\infty} \frac{(4A''')^{n}}{(n-1)n!}) \\
\leq &\frac{|\xi|}{N+n-2} (4A''' \log(N+n-2) + e^{4A'''}- 4A''' -1)\\
\leq &\frac{K \log (N+n-2)}{N+n-2}
\end{align*}
\end{proof}

\begin{proof}(Proof of theorem~\ref{FSA})
Here we assume the following lemmas which can be proved by standard techniques in complex analysis. Lemma~\ref{SE} gives an upper bound of the error terms appear in our estimation, while lemma~\ref{CML} provides a coordinate chart where explicit calculation can be done.

\begin{lemma}\label{SE}(Simple estimate)
Let $f: W\subset \CC \to \CC$ be a holomorphic function and $C$ be a contour in $W$. Let $M = sup_{z \in C} \Re (f(z)) < +\infty$. If there exists $N_{0}>0$ such that the integral $\int_{C} |\exp(N_{0}f(z))| dz$ is finite, then for $N>N_{0}$, we have
$$ | \int_{C} \exp(Nf(z)) dz | \leq C(f,N_{0}) \exp(NM),$$
where $C(f,N_{0}) = \exp(-N_{0}M) \int_{C} |\exp(N_{0}f(z))| dz $ is a constant depending on $N_{0}$ and $f$.
\end{lemma}

\begin{lemma}\label{CML}(Complex Morse lemma)
Let $f: W\subset \CC \to \CC$ be a holomorphic function and let $z_{0}$ be the only saddle point of $f$ in $W$. Further assume that the saddle point is non-degenerate. Then there exist a neighborhood $B(0,\delta) \subset \CC$ of $0$ with $\delta \in (0,1)$, a neighborhood $U(\delta) \subset W$ of $z$ with $B(z,\delta)\subset U(\delta)$ and a bijective holomorphic function $h:B(0,\delta) \to U(\delta)$ with $h(0)=z_{0}$ such that for any $w \in B(0,\delta)$,
\begin{equation}\label{pCML}
f(h(w)) = f(z_{0}) + \frac{1}{2} \frac{d^{2}f}{dz^{2}}(z_{0}) w^{2} \quad\quad \text{ and }\quad\quad  \frac{d h}{dw} (0) =1
\end{equation}
\end{lemma}

Now we outline the proof of the ordinary saddle point approximation and explicitly construct a constant coming from the term $O(1/N)$. To prove theorem~\ref{FSA}, it suffices to show that we can choose the constant to be independent on $y$ whenever $y$ is small. To do so, let us recall the statement of the saddle point approximation:

\begin{theorem}\label{SDA}(Saddle point approximation)
Let $f: W\subset \CC \to \CC$ be a holomorphic function and let $z_{0}$ be the only saddle point of $f$ in $W$. Further assume that the saddle point is non-degenerate and the maximum of the real part of $f$ attains at $z_{0}$. Let $C \subset W$ be a contour with finite length passing through the saddle point $z_{0}$. Assume that $|\arg( \sqrt{-\frac{d^{2}f}{dz^{2}}(z_{0})})| < \pi /4$. Then we have the following asymptotic formula:
\begin{equation}
I(N) = \int_{C} \exp(Nf(z)) dz =  \sqrt{\frac{2\pi}{N(-\frac{d^{2} f(z_{0})}{dz^{2}})}} \exp(N f(z_{0}))(1 + O(\frac{1}{N}))
\end{equation}
\end{theorem}

Given a function $f$ satisfying the properties stated in the theorem~\ref{SDA} , by lemma~\ref{CML} one can find neighborhoods $B(0,\delta)$ and $U(\delta)$ together with a bijective holomorphic function $h$ satisfying property (\ref{pCML}). Let $C_{0} = C \cap U$ and $C_{1} = C \backslash U$. We decompose the integral into two parts as follows:
$$I(N) = I_{0}(N)+ I_{1}(N) = \int_{C_{0}} \exp(Nf(z)) dz + \int_{C_{1}} \exp(Nf(z)) dz$$

Let $M = sup_{z \in C_{1}} \Re (f(z)) < +\infty$ and $l(C_{1})=\text{length of $C_{1}$}$. We have
\begin{equation}\label{EI1}
|I_{1}(N) |= | \int_{C_{1}} \exp(Nf(z)) dz | \leq l(C_{1}) \exp(NM)  ,
\end{equation}

Later we will show that this integral can be ignored when $N \to \infty$. So it suffices to consider the integral $I_{0}$. By change of variable formula one has
\begin{equation*}
I_{0}(N)=  \int_{C_{0}} \exp(Nf(z)) dz = \int_{h^{-1}(C_{0})} \exp(Nf(h(w))) \frac{dh}{dw}(w) dw
\end{equation*}

Recall that $z_{0} = h(0)$ and $f(h(w)) = f(z_{0}) + \frac{1}{2} \frac{d^{2}f}{dz^{2}}(z_{0}) w^{2}$. Consider the integration along the x-axis, i.e.
$$I'_{0}(N) = \int_{-\delta}^{\delta} \exp(N(f(z_{0}) + \frac{1}{2} \frac{d^{2}f}{dz^{2}}(z_{0}) x^{2})) \frac{dh}{dw}(x) dx$$

By analyticity of the integrand, the difference between $I_{0}(N)$ and $I'_{0}(N)$ can be expressed as
$$ E(N) =  \int_{\Gamma_{1}} \exp(Nf(z)) dz + \int_{\Gamma_{2}} \exp(Nf(z)) dz,$$
where $\Gamma_{i}$'s$\subset \partial B(0,\delta)$ are circular arcs connecting the endpoints of $h^{-1}(C_{0})$ and $[-\delta, \delta]$.

Let $K_{i} =  sup_{z \in B(0,\delta)\backslash B(0,\delta/2)} \Re (f(z)) < +\infty$. One can easily see that
\begin{align}\label{EE}
|\int_{\Gamma_{1}} \exp(Nf(z)) dz | &\leq 2\pi \delta \exp(NK_{1}) \text{ and } \\
|\int_{\Gamma_{1}} \exp(Nf(z)) dz | &\leq 2\pi \delta \exp(NK_{2}) ,
\end{align}

Furthermore, we extend the domain of integration to the whole real line. The error can be estimated by lemma~\ref{SE}. i.e. if  $Q = \sup_{x \in (-\infty, -\delta) \cup (\delta, \infty)} \Re (f(x)) $,
\begin{equation}\label{EG}
|G(N) |= | \int_{(-\infty, -\delta) \cup (\delta, \infty)} \exp(Nf(z)) dz | \leq P(f) \exp(NQ)  ,
\end{equation}
where $P(f) =  \exp(-Q) \int_{(-\infty, -\delta) \cup (\delta, \infty)} |\exp(f(z))| dz $.

Consider the Taylor's series expansion $ \frac{dh}{dw}(w)=1+ \sum_{n=1}^{\infty}a_{n}w^{n}$. We will first compute the contribution of the zero order term. The contribution of the higher order terms will be discussed later. By direct calculation we have
\begin{align*}
 &\int_{-\infty}^{\infty} \exp(N(f(z_{0}) + \frac{1}{2} \frac{d^{2}f}{dz^{2}}(z_{0}) x^{2})) dx \\
= & \exp(Nf(z_{0})) \int_{-\infty}^{\infty} \exp(\frac{1}{2} \frac{d^{2}f}{dz^{2}}(z_{0}) x^{2})) dx  \\
=& 2\exp(Nf(z_{0}))  \int_{0}^{\infty} \exp(-\frac{1}{2}N|\sqrt{-\frac{d^{2}f}{dz^{2}}(z_{0})}|^{2}x^{2} \exp(2i \arg \sqrt{-\frac{d^{2}f}{dz^{2}}})) dx
\end{align*}

Together with the assumption that $|\arg( \sqrt{-\frac{d^{2}f}{dz^{2}}(z_{0})})| < \pi /4$, the integral exists. Furthermore, by a change of variable we obtain
\begin{align}\label{MT}
&2\exp(Nf(z_{0}))  \int_{0}^{\infty} \exp(-\frac{1}{2}N|\sqrt{-\frac{d^{2}f}{dz^{2}}(z_{0})}|^{2}x^{2} \exp(2i \arg \sqrt{-\frac{d^{2}f}{dz^{2}}})) dx \notag \\
=& \sqrt{\frac{2\pi}{N(-\frac{d^{2} f(z_{0})}{dz^{2}})}} \exp(N f(z_{0}))
\end{align}

For the contribution of the higher order terms, note that when $n$ is odd, since the integrand is odd and the limit converges,
\begin{align*}
\int_{-\infty}^{\infty} \exp(N(f(z_{0}) + \frac{1}{2} \frac{d^{2}f}{dz^{2}}(z_{0}) x^{2}))x^{n} dx=0
\end{align*}

When $n$ is even, by integration by part, we have
\begin{align*}
&\int_{-\infty}^{\infty} \exp(N(f(z_{0}) + \frac{1}{2} \frac{d^{2}f}{dz^{2}}(z_{0}) x^{2}))x^{n} dx\\
=& \frac{1}{n+1}\int_{-\infty}^{\infty} \exp(N(f(z_{0}) + \frac{1}{2} \frac{d^{2}f}{dz^{2}}(z_{0}) x^{2})) dx^{n+1}\\
= &\frac{N\frac{d^{2}f}{dz^{2}}(z_{0})}{n+1}\int_{-\infty}^{\infty} \exp(N(f(z_{0}) + \frac{1}{2} \frac{d^{2}f}{dz^{2}}(z_{0}) x^{2})) x^{n+2}dx
\end{align*}

Iteratively, for positive integer $k$, the contribution of the degree $2k$ term is given by
\begin{align*}
&\int_{-\infty}^{\infty} \exp(N(f(z_{0}) + \frac{1}{2} \frac{d^{2}f}{dz^{2}}(z_{0}) x^{2})) x^{2k}dx \\
= &  (\frac{1}{N\frac{d^{2}f}{dz^{2}}(z_{0})})^{k}(\int_{-\infty}^{\infty} \exp(N(f(z_{0}) + \frac{1}{2} \frac{d^{2}f}{dz^{2}}(z_{0}) x^{2}))dx)
\end{align*}

In particular, the contribution of the second order term is given by
\begin{align}\label{SOT}
H(N) =&\int_{-\infty}^{\infty} \exp(N(f(z_{0}) + \frac{1}{2} \frac{d^{2}f}{dz^{2}}(z_{0}) x^{2})) x^{2}dx \notag \\
=& \frac{1}{N\frac{d^{2}f}{dz^{2}}(z_{0})}(\int_{-\infty}^{\infty} \exp(N(f(z_{0}) + \frac{1}{2} \frac{d^{2}f}{dz^{2}}(z_{0}) x^{2}))dx)
\end{align}

The sum of the contribution of the higher order terms is given by
\begin{align}\label{HOT}
&\sum_{k=2}^{\infty}a_{2k}(\frac{1}{N\frac{d^{2}f}{dz^{2}}(z_{0})})^{k}(\int_{-\infty}^{\infty} \exp(N(f(z_{0}) + \frac{1}{2} \frac{d^{2}f}{dz^{2}}(z_{0}) x^{2}))dx) \notag \\
= & \sum_{k=2}^{\infty}a_{2k}(\frac{1}{N\frac{d^{2}f}{dz^{2}}(z_{0})})^{k-1} H(N)
\end{align}

By comparing equation (\ref{SOT}) and (\ref{HOT}), one can see that the contribution of the higher order terms can be ignored compared with that of second order term.

Furthermore, for $I_{1}(N)$, $E(N)$ and $G(N)$, from inequalities (\ref{EI1}), (\ref{EE}) and (\ref{EG}), they grow in $\exp (N \times \text{constant})$, where the constant is strictly less than $\Re (f(z_{0}))$. As a result, these error terms decay exponentially compared with equation (\ref{MT}).

To conclude, we may take the constant appear in $O(1/N)$ to be $ 2/\frac{d^{2}f}{dz^{2}}(z_{0})$ such that whenever $N$ is large, we have
\begin{align*}
|\int_{C} \exp(Nf(z)) dz / \sqrt{\frac{2\pi}{N(-\frac{d^{2} f(z_{0})}{dz^{2}})}} \exp(N f(z_{0})) - 1| \leq \frac{2/\frac{d^{2}f}{dz^{2}}(z_{0})}{N}
\end{align*}

Finally we can prove theorem~\ref{FSA}. That means we have to control the error terms uniformly on $y$. To do so, for each $y \in [0,1]$ we apply lemma~\ref{CML} for $\Phi_{y}(z)$ to find $U_{y}(\delta_{y})$ containing the saddle point $z_{y}$. From theorem 1.1 in \cite{CHP03}, one can check that the size of the neighborhood in theorem~\ref{CML} has a lower bound which depends continuously on the function $f$. Therefore in our situation we can find a $\delta>0$ such that $U_{y_{N}}(\delta) \subset U_{y_{N}}(\delta_{N})$ for all sufficiently large $N$. From this we can find a good control of the supremum of our functions $\Phi_{y_{N}}$ outside $U_{y_{N}}(\delta)$ as follows.

Apply the same argument to each function $\Phi_{y_{N}}$ as in the proof of the ordinary saddle point approximation and denote the constants appeared in the estimation (\ref{EI1}), (\ref{EE}) and (\ref{EG}) to be $M(\Phi_{y_{N}})$, $K_{1}(\Phi_{y_{N}})$, $K_{2}(\Phi_{y_{N}})$ and $Q(\Phi_{y_{N}})$ respectively. By the continuity of the function $h(y,t)=\Phi_{y}(C(y,t)): [0,1]^{2} \to \CC$, for $N$ sufficiently large we have
\begin{align*}
&|M(\Phi_{y_{N}}) - \Phi_{y_{N}}(z_{y_{N}})| \\
\geq & -|M(\Phi_{y_{N}}) - M(\Phi_{0})| + |M(\Phi_{0}) - \Phi_{0}(z_{0})| - |\Phi(z_{0}) - \Phi_{y_{N}}(z_{y_{N}})| \\
\geq & |M(\Phi_{0}) - \Phi_{0}(z_{0})|/2 >0
\end{align*}

Moreover, by our assumption the length of the contours $C_{y}$ are uniformly bounded by a constant $L$. This provides a uniform way for exponential decay.

Similar arguments can be applied to $K_{1}(\Phi_{y_{N}}), K_{2}(\Phi_{y_{N}})$ and $Q(\Phi_{y_{N}})$. Thus the constants can be chosen to be independent on $N$ whenever $N$ is large.

Moreover, the coefficient of $H(N)$ appeared in equation (\ref{HOT}) depend smoothly on the function $f$. Under the assumption that $\Phi_{y_{N}} \xrightarrow{N \to \infty} \Phi_{0}$ the constant can be chosen to be independent on $N$. Together with the fact that $\ds \frac{d^{2}\Phi_{y_{N}}}{dz^{2}}(z_{y_{N}}) \to \frac{d^{2}\Phi_{0}}{dz^{2}}(z_{0})$ as $N \to \infty$, the error term can be chosen uniformly on $N$. This completes the proof of theorem~\ref{FSA}.
\end{proof}

\begin{proof}(Proof of theorem~\ref{saddlepointapp})
To prove theorem~\ref{saddlepointapp}, it suffices to show that the conditions in theorem~\ref{FSA} are satisfied in our situation.\\

First of all we show the existence of such paths when $N$ is sufficiently large. We are going to construct the contour using the same idea as in the proof of lemma 3.4 of \cite{HM13}. To do so, we only need to check that the conditions in the construction are also satisfied in our case.\\

Let $q_{N}(t)=z^{(n)}_{N}t$ for $0<t<\operatorname{Re}(1/z^{(n)}_{N})$. Since $\ds \lim_{N \to \infty}z^{(n)}_{N}=z^{(2)}<1$ (see the proof of lemma 3.4 of \cite{HM13}), $\operatorname{Re}(1/z^{(n)}_{N})>1$ for sufficiently large $N$. Also, since $d^{2}\Phi^{(2)}(z^{(2)})/dz^{2} \neq 0$ and $\Phi^{(n)}_{N} \to \Phi^{(2)}$ as $N$ goes to infinity, we have $d^{2}\Phi^{(n)}_{N}(z^{(n)}_{N})/d z^{2} \neq 0$  for sufficiently large $N$. By definition we have $d \Phi^{(n)}_{M}(z^{(n)}_{M}) /dz =0$. This implies $\operatorname{Re}\Phi^{(n)}_{N}(q_{N}(1))=0$ for any $N$. Since $\max \{\operatorname{Re}\Phi^{(2)}(z)\}$ takes place at $z=z^{(2)}$, we must have $\max \{\operatorname{Re}\Phi^{(n)}_{N}(z)\} = \operatorname{Re}\Phi^{(n)}_{N}(z^{(n)}_{N})$ along the line $q_{N}(t)$. \\

Moreover, from the proof of lemma 3.4 of \cite{HM13} that the difference between the argument of $z^{(2)}$ and $1/\sqrt{-d^{2}\Phi^{(2)}(z^{(2)})/dz^{2}}$ is strictly smaller than $\pi/4$. Hence the difference between the argument of $z^{(n)}_{N}$ and $1/\sqrt{-d^{2}\Phi^{(n)}_{N}(z^{(n)}_{N})/dz^{2}}$ is also strictly smaller than $\pi/4$ for large $N$. As a result the same construction of the path $Q$ in the proof of lemma 3.4 of \cite{HM13} still applies. \\

Finally we connect $z^{(n)}_{N}(\operatorname{Re}1/z^{(n)}_{N})$ and $1$ by a line segment $L$. Since from the proof of lemma 3.4 in \cite{HM13} that $\operatorname{Re}\Phi^{(2)}(z)<0$ on the segment connecting $2\pi i/\xi$ and $1$, $\operatorname{Re}\Phi^{(n)}_{N}(\omega)\leq 0$ on the segment $L$ for large $N$. This finishes the construction of the paths. We will denote the contours by $Q_{N}$.\\

Theorem~\ref{saddlepointapp} follows from direct application of theorem~\ref{FSA} with the data
\begin{align*}
\Phi^{(n)}_{y}(z)&=\frac{1}{\xi}(\operatorname{Li}_{2}(e^{u-(z+y\xi}) + \operatorname{Li}_{2}(e^{z\xi}) - \operatorname{Li}_{2}(e^{u+z\xi}) - \operatorname{Li}_{2}(e^{(z+y\xi}))-uz \\
y_{N}&=\frac{n-2}{N} \text{ and } C_{y_{N}}=Q_{N}
\end{align*}
\end{proof}

\begin{proof}(Proof of lemma~\ref{Sratio})
\begin{align*}
&\frac{S_{\gamma}(-\pi - iu + \gamma)}{S_{\gamma}(\pi - iu - (2n-3)\gamma)} \\
= &\exp \left(\frac{1}{4} \int_{C_{R}} \frac{e^{-iut}}{\sinh (\pi t) \sinh (\gamma t)}(e^{-\pi t + \gamma t} - e^{\pi t - (2n-3)\gamma t}) \frac{dt}{t}\right) \\
= &\exp \left(\frac{1}{4} \int_{C_{R}} \frac{e^{-iut}e^{-(n-2)\gamma t}}{\sinh (\pi t) \sinh (\gamma t)}(e^{-\pi t + (n-1)\gamma t} - e^{\pi t - (n-1)\gamma t}) \frac{dt}{t}\right) \\
= &\exp \left(\frac{1}{2} \int_{C_{R}} e^{-iut}e^{-(n-2)\gamma t} \frac{ \sinh(-\pi t + (n-1)\gamma t)}{\sinh (\pi t) \sinh (\gamma t)} \frac{dt}{t}\right) \\
=  &\exp \left( \frac{1}{2} \int_{C_{R}} \frac{e^{-iut}\coth(\pi t)}{t} (e^{-(n-2)\gamma t} \frac{\sinh(n-1)\gamma t}{\sinh (\gamma t)})  -  \right. \\
&\quad \quad \quad \quad \quad \left. \frac{e^{-iut}\coth(\gamma t)}{t} (e^{-(n-2)\gamma t}\frac{\cosh(n-1)\gamma t}{\cosh (\gamma t)}) dt \right)
\end{align*}

Furthermore, one can easily verify the following formulas: (the proofs will be given later)
\begin{equation}
\begin{aligned}
e^{-(n-2)A}\frac{\sinh((n-1)A)}{\sinh(A)} &= \sum_{k=0}^{n-2}e^{-2kA}\\ \label{TL1}
e^{-(n-2)A}\frac{\cosh((n-1)A)}{\cosh(A)} &= \sum_{k=0}^{n-2}(-1)^{k}e^{-2kA}
\end{aligned}
\end{equation}

From these formulas we can see that
\begin{align*}
&\frac{S_{\gamma}(-\pi - iu + \gamma)}{S_{\gamma}(\pi - iu - (2n-3)\gamma)}\\
= &\exp \left(\frac{1}{4} \int_{C_{R}} \frac{e^{-iut}}{\sinh (\pi t) \sinh (\gamma t)}(e^{-\pi t + \gamma t} - e^{\pi t - (2n-3)\gamma t}) \frac{dt}{t} \right) \\
=  &\exp \left(\frac{1}{2} \int_{C_{R}} \frac{e^{-iut}\coth(\pi t)}{t} (\sum_{k=0}^{n-2}e^{-2k\gamma t}) - \frac{e^{-iut}\coth(\gamma t)}{t} (\sum_{k=0}^{n-2}(-1)^{k}e^{-2k\gamma t}) dt \right)
\end{align*}

Now we modify the proof in \cite{HM13}. For $r>0$, let $U_{i}, i = 1,2,3$ be the segments defined by $r \xrightarrow{U_{1}} r - r' i \xrightarrow{U_{2}}  -r  - r' i \xrightarrow{U_{3}}  -r$ with $r' = \frac{3\pi}{u}r$. Since the zeros of $\sinh(\pi t)$ and $\sinh(\gamma t)$  are discrete, for genreric $r'$, $U_{2}$ does not pass through those singular points.

Now we want to show that for $i=1,2,3$,
\begin{align*}
\lim_{r \to \infty} \int_{U_{i}} \frac{e^{-iut}}{\sinh (\pi t) \sinh (\gamma t)}(e^{-\pi t + \gamma t} - e^{\pi t - (2n-3)\gamma t}) \frac{dt}{t} = 0
\end{align*}

We will show the convergence on (i) $U_{1}$, (ii) $U_{3}$, (iii) $U_{2}$.\\

First of all we choose $r>0$ satisfying
\begin{align*}
r = \frac{(2l+1)\pi}{4\pi^{2}/(N+n-2)u + u/(N+n-2)} \text{ for $l \in \mathbb{N}$.}
\end{align*}
The choice of $r$ helps us to avoid the pole of $\sinh(\gamma t)$ and get a good estimation of the integrals.
More precisely, for $s\in [0,r']$ we consider the functions
$$p(s) = |1-e^{-2\pi(r-si)}|,\quad q(s) = |e^{2\pi(r-si)} - 1 | \text{ and } g(s) =  | e^{-2\gamma(r-si)} -1|$$
In the above $g(s)$ is the distance between $e^{-2\gamma(r-si)}$ and $1$. These functions correspond to the terms appear in the integrals as shown later. When $r$ is large,
\begin{align*}
p(s)&= |1-e^{-2\pi(r-si)}| \geq 1 - e^{-2\pi r} \geq 1/2 ;\\
q(s)&= |e^{2\pi(r-si)} - 1 | \geq e^{2\pi r} -1 \geq 1
\end{align*}
Also, one can check that $$g(s) =  | e^{-2\gamma(r-si)} -1| = | e^{R(s)}e^{i\theta(s)} - 1| ,$$
where $R(s)= \frac{us}{N+n-2} - \frac{2\pi r}{N+n-2}$ and $\theta(s)=\frac{ur}{N+n-2} + \frac{2\pi s}{N+n-2} $. Moreover, due to the choice of $r$,
\begin{itemize}
\item when $s = \frac{2\pi}{u}r$, we have $R(s)=0$, $\theta(s)=(2l+1)\pi$;
\item when $s = \frac{2\pi}{u}r - \frac{N+n-2}{4}$, we have $R(s)=-\frac{u}{4}$, $\theta(s)=(2l+1)\pi - \frac{\pi}{2}$;
\item when $s = \frac{2\pi}{u}r + \frac{N+n-2}{4}$, we have $R(s)=\frac{u}{4}$, $\theta(s)=(2l+1)\pi + \frac{\pi}{2}$.
\end{itemize}

Since $R(s)$ and $\theta(s)$ are strictly increasing in $s$ and $g(s)$ is the distance between $ e^{R(s)}e^{i\theta(s)}$ and $1$,
\begin{itemize}
\item for $0 \leq s \leq \frac{2\pi}{u}r - \frac{N+n-2}{4}$, $\ds g(s) \geq \min_{|z|\leq e^{-u/4}} |z -1| = 1 - e^{-u/4} $
\item for $\frac{2\pi}{u}r - \frac{N+n-2}{4} \leq s \leq \frac{2\pi}{u}r + \frac{N+n-2}{4}$, since $\theta(s) \in [(2l+1)\pi - \frac{\pi}{2}, (2l+1)\pi + \frac{\pi}{2}]$, we must have $g(s) \geq 1$.\\
\item for $\frac{2\pi}{u}r + \frac{N+n-2}{4} \leq s \leq  \frac{3\pi}{u}r$, $\ds g(s) \geq  \min_{|z|\geq e^{u/4}} |z -1| = e^{u/4} -1$.
\end{itemize}

To conclude, we can find positive constants $M_{1}$, $M_{2}$ and $M_{3}$ independent on $r$ such that
$$ \frac{1}{p(s)} \leq M_{1}, \quad \frac{1}{q(s)} \leq M_{2}, \quad \frac{1}{g(s)} \leq M_{3} $$

Now we can get a good control of the integrals.\\
(i) On $U_{1}$,
\begin{align*}
&|\int_{U_{1}} \frac{e^{-iut}}{\sinh (\pi t) \sinh (\gamma t)}(e^{-\pi t + \gamma t}) \frac{dt}{t})|\\
\leq& 4 \int_{0}^{r'} |\frac{e^{-iu(r-si)}}{r-si}| |\frac{e^{(-\pi + \gamma) (r-si)}}{\sinh(\pi (r-si))\sinh(\gamma (r-si))}| ds \\
\leq& \frac{4M_{2}M_{3}}{r} \int_{0}^{r'} e^{-us} ds\\
=& \frac{4M_{2}M_{3}}{ur}(1-e^{-ur'}) \xrightarrow{r \to \infty} 0.
\end{align*}

Similarly,
\begin{align*}
&|\int_{U_{1}} \frac{e^{-iut}}{\sinh (\pi t) \sinh (\gamma t)}(e^{\pi t - (2n-3)\gamma t}) \frac{dt}{t})|\\
\leq& 4\int_{0}^{r'} |\frac{e^{-iu(r-si)}}{r-si}| |\frac{e^{(\pi - (2n-3)\gamma) (r-si)}}{\sinh(\pi (r-si))\sinh(\gamma (r-si))}| ds \\
\leq & \frac{4}{r} \int_{0}^{r'} e^{-us} |e^{-(2n-2)\gamma(r-si)}| \frac{1}{p(s)}\frac{1}{g(s)} ds\\
\leq & \frac{4M_{1}M_{3}}{r} \int_{0}^{r'} e^{(-1 + \frac{2n-2}{N+n-2})us - \frac{(2n-2)\pi r}{N+n-2}} ds
\end{align*}

Hence
\begin{align*}
&|\int_{U_{1}} \frac{e^{-iut}}{\sinh (\pi t) \sinh (\gamma t)}(e^{\pi t - (2n-3)\gamma t}) \frac{dt}{t})| \\
\leq & \frac{4M_{1}M_{3}}{r} \int_{0}^{r'} e^{(-1 + \frac{2n-2}{N+n-2})us - \frac{(2n-2)\pi r}{N+n-2}} ds \\
\leq & \frac{4M_{1}M_{3}}{ur}(e^{(-1 + \frac{2n-2}{N+n-2})ur' - \frac{(2n-2)\pi r}{N+n-2}} - e^{-\frac{(2n-2)\pi r}{N+n-2}})  \xrightarrow{r \to \infty} 0,
\end{align*}

(ii) On $U_{3}$,
\begin{align*}
&|\int_{U_{3}} \frac{e^{-iut}}{\sinh (\pi t) \sinh (\gamma t)}(e^{-\pi t + \gamma t}) \frac{dt}{t})|\\
\leq& 4 \int_{0}^{r'} |\frac{e^{-iu(-r-si)}}{-r-si}| |\frac{e^{(-\pi + \gamma)(-r-si)}}{(e^{\pi(-r-si)}-e^{-\pi(-r-si)})(e^{\gamma(-r-si)}-e^{-\gamma(-r-si)})}| dt \\
\leq& \frac{4}{r} \int_{0}^{r'} e^{-us} \frac{1}{|e^{2\pi(-r-si)} + e^{-2\gamma{(-r-si)}} - 1 - e^{2(\pi-\gamma)(-r-si)}|} dt
\end{align*}
Note that the modulus of the terms in the denominator are 
$$e^{-2\pi r}, e^{\frac{2\pi r}{N+n-2}+\frac{us}{N+n-2}}, 1 \text{ and } e^{-2\pi r + \frac{2\pi r}{N+n-2}+\frac{us}{N+n-2}}$$
respectively. For large r, the dominant term is $e^{\frac{2\pi r}{N+n-2}+\frac{us}{2(N+n-2)}} \xrightarrow{r \to \infty} \infty$. This show that the denominator is bounded below. So we can find some constant $M_{4}$ such that
\begin{align*}
|\int_{U_{3}} \frac{e^{-iut}}{\sinh (\pi t) \sinh (\gamma t)}(e^{-\pi t + \gamma t}) \frac{dt}{t})|
\leq\frac{M_{4}}{r} \int_{0}^{r'} e^{-us} dt
\leq \frac{M_{4}}{ur}(1-e^{-ur'}) \xrightarrow{r \to \infty} 0.
\end{align*}

Similarly,
\begin{align*}
&|\int_{U_{3}} \frac{e^{-iut}}{\sinh (\pi t) \sinh (\gamma t)}(e^{\pi t - (2n-3)\gamma t}) \frac{dt}{t})|\\
\leq& 4 \int_{0}^{r'} |\frac{e^{-iu(-r-si)}}{-r-si}| |\frac{e^{(\pi - (2n-3) \gamma)(-r-si)}}{(e^{\pi(-r-si)}-e^{-\pi(-r-si)})(e^{\gamma(-r-si)}-e^{-\gamma(-r-si)})}| dt \\
\leq & \frac{4}{r} \int_{0}^{r'} e^{-us} \frac{1}{|e^{q_{1}} + e^{q_{2}} -e^{q_{3}} -e^{q_{4}}|}
\end{align*}
where $q_{1}$, $q_{2}$, $q_{3}$ and $q_{4}$ are given by
\begin{align*}
q_{1}&=(2n-2)\gamma(-r-si), \quad \quad \quad \quad\quad \quad  q_{2}=(-2\pi+(2n-4)\gamma)(-r-si),\\
q_{3}&=(-2\pi+(2n-2)\gamma)(-r-si),\quad \quad  q_{4}=((2n-4)\gamma)(-r-si)
\end{align*}
Note that the modulus of the terms in the denominator are 
$$e^{-\frac{2(n-1)\pi r}{N+n-2}- \frac{(n-1)us}{N+n-2}}, e^{2\pi r - \frac{2(n-2)\pi r}{N+n-2}- \frac{(n-2)us}{N+n-2}}, e^{2\pi r - \frac{2(n-1)\pi r}{N+n-2}- \frac{(n-1)us}{N+n-2}} \text{ and } e^{-\frac{2(n-2)\pi r}{N+n-2}- \frac{(n-2)us}{N+n-2}}$$ respectively. For large r, the dominant term is $e^{2\pi r - \frac{2(n-2)\pi r}{N+n-2}- \frac{(n-2)us}{2(N+n-2)}} \to \infty$. This show that the denominator is bounded below. Again we can find some constant $M_{5}$ such that
\begin{align*}
|\int_{U_{3}} \frac{e^{-iut}}{\sinh (\pi t) \sinh (\gamma t)}(e^{\pi t - (2n-3)\gamma t}) \frac{dt}{t})|
\leq \frac{M_{5}}{r} \int_{0}^{r'} e^{-us} dt 
\leq \frac{M_{5}}{ur}(1-e^{-ur'}) \xrightarrow{r \to \infty} 0.
\end{align*}

(iii) On $U_{2}$, we consider the expression
\begin{align*}
&\frac{S_{\gamma}(-\pi + iu + \gamma)}{S_{\gamma}(\pi - iu - (2n-3)\gamma)} \\
=  &\exp (\frac{1}{2} \int_{C_{R}} \frac{e^{-iut}\coth(\pi t)}{t} (\sum_{k=0}^{n-2}e^{-2k\gamma t}) - \frac{e^{-iut}\coth(\gamma t)}{t} (\sum_{k=0}^{n-2}(-1)^{k}e^{-2k\gamma t}) dt)
\end{align*}
Note that for $t=s - r'i$, $s \in [-r,r]$,
\begin{align*}
|e^{-2k\gamma t}|=e^{k(\frac{-2\pi s}{N+n-2} + \frac{r'u}{(N+n-2)})} \leq e^{\frac{k(2\pi r + r'u)}{(N+n-2)}} \leq e^{\frac{k(4\pi+u) r'}{(N+n-2)}}
\end{align*}
Write $\kappa=\alpha - \beta i$, where $\kappa = \pi$ or $\gamma$,
\begin{align*}
 |\int_{U_{2}} \frac{e^{-iut}\coth(\kappa t)}{t} (\sum_{k=0}^{n-2}e^{-2k\gamma t}) dt |
\leq &\int_{U_{2}} |\frac{e^{-iut}\coth(\kappa t)}{t}| (\sum_{k=0}^{n-2}|e^{-2k\gamma t}|)dt \\
\leq &\sum_{k=0}^{n-2}\frac{e^{-ur'(1-\frac{k(4\pi+u) }{u(N+n-2)})}}{r'} \int_{U_{2}} |\coth(\kappa t)| dt
\end{align*}

By the similar trick in \cite{HM13}, put $\delta = \max_{-1 \leq s \leq 1} |\coth(\kappa s)| > 0$. This helps us to get away from the singularity of $\coth(s\pi)$ in the proof shown below. Now we have
\begin{align*}
&\int_{U_{2}}|\coth(\kappa t)| dt \\
=& \int_{-r}^{r}|\coth(s\alpha - r'\beta - (s\beta + \alpha r')i)| ds \\
\leq & 2\delta + \int_{-r}^{-1} | \frac{e^{s\alpha - r'\beta - (s\beta + \alpha r')i} + e^{-(s\alpha - r'\beta - (s\beta + \alpha r')i)}}{e^{s\alpha - r'\beta - (s\beta + \alpha r')i} - e^{-(s\alpha - r'\beta - (s\beta + \alpha r')i)}} | ds \\
&+ \int_{1}^{r} | \frac{e^{s\alpha - r'\beta - (s\beta + \alpha r')i} + e^{-(s\alpha - r'\beta - (s\beta + \alpha r')i)}}{e^{s\alpha - r'\beta - (s\beta + \alpha r')i} - e^{-(s\alpha - r'\beta - (s\beta + \alpha r')i)}} | ds \\
\leq & 2\delta + \int_{-r}^{-1}  \frac{|e^{s\alpha - r'\beta - (s\beta + \alpha r')i}| + |e^{-(s\alpha - r'\beta - (s\beta + \alpha r')i)}|}{|e^{s\alpha - r'\beta - (s\beta + \alpha r')i}| - |e^{-(s\alpha - r'\beta - (s\beta + \alpha r')i)}|}  ds \\
& + \int_{1}^{r}  \frac{|e^{s\alpha - r'\beta - (s\beta + \alpha r')i}| + |e^{-(s\alpha - r'\beta - (s\beta + \alpha r')i)}|}{|e^{s\alpha - r'\beta - (s\beta + \alpha r')i}| - |e^{-(s\alpha - r'\beta - (s\beta + \alpha r')i)}|}  ds\\
= &2\delta +  \int_{-r}^{-1} \frac{|e^{s\alpha-r'\beta}| + |e^{-(s\alpha - r'\beta)}|}{|e^{s\alpha - r'\beta }| - |e^{-(s\alpha - r'\beta)}|}
+ \int_{1}^{r}  \frac{|e^{s\alpha-r'\beta}| + |e^{-(s\alpha - r'\beta)}|}{|e^{s\alpha - r'\beta }| - |e^{-(s\alpha - r'\beta)}|}  ds \\
\leq & 2\delta + \int_{1}^{r} \coth(s\alpha - r'\beta) ds +\int_{-r}^{-1} \coth(s\alpha - r'\beta) ds \\
= & 2\delta + \frac{\log(\sinh(\alpha r - r'\beta)) - \log(\sinh(\alpha - r'\beta))}{\alpha} \\
&+\frac{\log(\sinh(-\alpha - r'\beta)) - \log(\sinh(-\alpha r - r'\beta))}{\alpha}
\end{align*}

Hence
\begin{align*}
& |\int_{U_{2}} \frac{e^{-iut}\coth(\kappa t)}{t} (\sum_{k=0}^{n-2}e^{-2k\gamma t}) dt |\\
\leq &\sum_{k=0}^{n-2}\frac{e^{-ur'(1-\frac{k(4\pi+u) }{u(N+n-2)})}}{r'} [2\delta + \frac{\log(\sinh(\alpha r - r'\beta)) - \log(\sinh(\alpha - r'\beta))}{\alpha} \\
&+\frac{\log(\sinh(-\alpha - r'\beta)) - \log(\sinh(-\alpha r - r'\beta))}{\alpha}] \xrightarrow{r \to \infty} 0
\end{align*}

Let $C_{r}=[-r, -R] \cup \Omega_{R} \cup [R,r]$. Denote $U_{1}\cup U_{2} \cup U_{3}$ by $U_{123}$. By (i)-(iii) we get

\begin{align*}
&\int_{C_R} \frac{e^{-iut}}{\sinh (\pi t) \sinh (\gamma t)}(e^{-\pi t + \gamma t} - e^{\pi t - (2n-3)\gamma t}) \frac{dt}{t}) \\
=&\lim_{r \to \infty} \int_{C_{r}} (\frac{e^{-iut}}{\sinh (\pi t) \sinh (\gamma t)}(e^{-\pi t + \gamma t} - e^{\pi t - (2n-3)\gamma t}) \frac{dt}{t})) \\
=&\lim_{r \to \infty} \int_{C_{r}} (\frac{e^{-iut}\coth(\pi t)}{t} (\sum_{k=0}^{n-2}e^{-2k\gamma t}) - \frac{e^{-iut}\coth(\gamma t)}{t} (\sum_{k=0}^{n-2}(-1)^{k}e^{-2k\gamma t})) dt \\
=&\lim_{r \to \infty} [\int_{U_{123} } (\frac{e^{-iut}\coth(\pi t)}{t} (\sum_{k=0}^{n-2}e^{-2k\gamma t}) - \frac{e^{-iut}\coth(\gamma t)}{t} (\sum_{k=0}^{n-2}(-1)^{k}e^{-2k\gamma t}) )dt \\
&+2\pi i\text{Res}(\frac{e^{-iut}\coth(\pi t)}{t} (\sum_{k=0}^{n-2}e^{-2k\gamma t}), t=li) \\
& - 2\pi i\text{Res}(\frac{e^{-iut}\coth(\gamma t)}{t} (\sum_{k=0}^{n-2}(-1)^{k}e^{-2k\gamma t}), t=\frac{l\pi i}{\gamma}) ]\\
=&\lim_{r \to \infty} \int_{U_{123}} \frac{e^{-iut}}{\sinh (\pi t) \sinh (\gamma t)}(e^{-\pi t + \gamma t} - e^{\pi t - (2n-3)\gamma t}) \frac{dt}{t} \\
&+2\pi i\text{Res}(\frac{e^{-iut}\coth(\pi t)}{t} (\sum_{k=0}^{n-2}e^{-2k\gamma t}), t=li) \\
&- 2\pi i\text{Res}(\frac{e^{-iut}\coth(\gamma t)}{t} (\sum_{k=0}^{n-2}(-1)^{k}e^{-2k\gamma t}),t=\frac{l\pi i}{\gamma})\\
=&\lim_{r \to \infty}2\pi i\text{Res}(\frac{e^{-iut}\coth(\pi t)}{t} (\sum_{k=0}^{n-2}e^{-2k\gamma t}), t=li) \\
&- 2\pi i\text{Res}(\frac{e^{-iut}\coth(\gamma t)}{t} (\sum_{k=0}^{n-2}(-1)^{k}e^{-2k\gamma t}),t=\frac{l\pi i}{\gamma})
\end{align*}

Here ``Res'' means all the residue inside the contour $C_{r} \cup U_{1} \cup U_{2}\cup U_{3}$ as $r$ goes to infinity. So we have

\begin{align*}
&\int_{C_{R}} \frac{e^{-iut}e^{-(n-2)\gamma t}}{t} \coth(\pi t) \frac{\sinh(n-1)\gamma t}{\sinh (\gamma t)} \\
=& 2\pi i \sum_{l=1}^{\infty} \text{Res} (\frac{e^{-iut}e^{-(n-2)\gamma t}}{t} \coth(\pi t) \frac{\sinh(n-1)\gamma t}{\sinh (\gamma t)} ; t = li) \\
= &2\pi i \sum_{l=1}^{\infty}\frac{e^{ul}}{l\pi i}\sum_{k=0}^{n-2}e^{-2k\gamma l i} \\
= &2\pi i \sum_{k=0}^{n-2}\sum_{l=1}^{\infty}\frac{e^{l(u-2k\gamma i)}}{l\pi i} \\
= &-2\sum_{k=0}^{n-2} \log (1-e^{u-2k\gamma i})
\end{align*}

and
\begin{align*}
&\int_{C_{R}} \frac{e^{-iut}e^{-(n-2)\gamma t}}{t} \coth(\gamma t)\frac{\cosh(n-1)\gamma t}{\cosh (\gamma t)} \\
= &2\pi i \sum_{l=1}^{\infty} \text{Res} (\frac{e^{-iut}e^{-(n-2)\gamma t}}{t} \coth(\gamma t) \frac{\cosh(n-1)\gamma t}{\cosh (\gamma t)} ; t =\frac{l\pi i}{\gamma}) \\
= &2\pi i \sum_{l=1}^{\infty}\frac{e^{ul\pi/\gamma}}{l\pi i}\sum_{k=0}^{n-2}(-1)^{k}e^{-2k l\pi i} \\
= & 2\pi i \sum_{l=1}^{\infty}\frac{e^{ul\pi/\gamma}}{l\pi i} \\
=& -2 \log ( 1 - e^{u\pi / \gamma} )
\end{align*}

Overall we have
\begin{align*}
\frac{S_{\gamma}(-\pi-iu+\gamma)}{S_{\gamma}(\pi-iu-(2n-3)\gamma)}
&=
\frac{ e^{u\pi / \gamma}-1 }{ \prod_{k=0}^{n-2} (e^{u-2k\gamma i}-1) } \stackrel[N \to \infty]{\sim}{ }
\frac{e^{2\pi i u N / \xi}}{(e^{u}-1)^{n-1}}
\end{align*}

\end{proof}

\begin{proof}(Proof of equations (\ref{TL1}))
\begin{align*}
e^{-(n-2)A}\frac{\sinh((n-1)A)}{\sinh(A)}
&= e^{-(n-2)A}\frac{e^{(n-1)A} - e^{-(n-1)A}}{e^{A}-e^{-A}} \\
&= \frac{1-e^{-2(n-1)A}}{1-e^{-2A}} \\
&= \frac{1-(e^{-2A})^{n-1}}{1-e^{-2A}} \\
&= \sum_{k=0}^{n-2}e^{-2kA}
\end{align*}
The second equality requires the condition that $n$ is even. The proof is similar as above so we omit it.
\end{proof}g

\begin{proof}(Proof of proposition~\ref{positive})
From Lemma 3.5 in \cite{HM13} we know that 
$$\operatorname{Re}\Phi^{(2)}(\omega_{0})>0 \text{ for } 0<u<\log((3+\sqrt{5})/2).$$
 Since $\Phi^{(n)}_{N}(\omega^{(n)}_{N}) \to \Phi(\omega)$ as $N \to \infty$, we get the result.
\end{proof}

\begin{proof} (Proof of lemma~\ref{diff1})
Recall that
\begin{align*}
\Phi ^{(2)}(z) &= \frac{1}{\xi}(\operatorname{Li}_{2}(e^{u-\xi z}) - \operatorname{Li}_{2}(e^{u+\xi z})) - uz \\
\Phi^{(n)}_{N}(z) &=\frac{1}{\xi}(\operatorname{Li}_{2}(e^{u-(z+\frac{n-2}{N})\xi}) + \operatorname{Li}_{2}(e^{z\xi}) - \operatorname{Li}_{2}(e^{u+z\xi}) - \operatorname{Li}_{2}(e^{(z+\frac{n-2}{N})\xi}))-uz
\end{align*}
So
\begin{align*}
\Phi^{(n)}_{N}(z) - \Phi ^{(2)}(z) =  \frac{1}{\xi} ( \operatorname{Li}_{2}(e^{u-(z+\frac{n-2}{N})\xi}) - \operatorname{Li}_{2}(e^{u-z\xi }) + \operatorname{Li}_{2}(e^{z\xi}) -  \operatorname{Li}_{2}(e^{(z+\frac{n-2}{N})\xi}) )
\end{align*}
Put $\ds y = \frac{n-2}{N}$, we have
\begin{align*}
\Phi^{(n)}_{N}(z) - \Phi ^{(2)}(z) =  \frac{1}{\xi} ( \operatorname{Li}_{2}(e^{u-(z+y)\xi}) - \operatorname{Li}_{2}(e^{u-z\xi }) + \operatorname{Li}_{2}(e^{z\xi}) -  \operatorname{Li}_{2}(e^{(z+y)\xi}) )
\end{align*}
\begin{align*}
&\lim_{N \to \infty} N(\Phi_{N}^{(n)}(z) - \Phi^{(2)}(z)) \\
=& \frac{n-2}{\xi} \lim_{y \to 0} \frac{ \operatorname{Li}_{2}(e^{u-(z+y)\xi}) - \operatorname{Li}_{2}(e^{u-z\xi }) + \operatorname{Li}_{2}(e^{z\xi}) -  \operatorname{Li}_{2}(e^{(z+y)\xi})}{y} \\
=& \frac{n-2}{\xi} \lim_{y \to 0} \frac{d}{dy} ( \operatorname{Li}_{2}(e^{u-(z+y)\xi}) - \operatorname{Li}_{2}(e^{u-z\xi }) + \operatorname{Li}_{2}(e^{z\xi}) -  \operatorname{Li}_{2}(e^{(z+y)\xi})) \\&(\text{by L'Hospital's rule}) \\
=&  \frac{n-2}{\xi} \lim_{y \to 0} ( - \log(1-e^{u-(z+y)\xi}) (- \xi) - (- \log(1-e^{(z+y)\xi}))(\xi)) \\
=& (n-2)\log((1-e^{u-z\xi})(1-e^{z\xi}))
\end{align*}
\end{proof}

\begin{proof}(Proof of lemma~\ref{diff2})
To remove the $N$ dependence of $z^{(n)}_{N+n-2}$, recall that from (\ref{cubicpoly})
\begin{equation*}
ab^{2}(\omega^{(n)}_{N+n-2})^{3} - (b^{2}+a^{2}b)(\omega^{(n)}_{N+n-2})^{2} + (a^{2}+b)(\omega^{(n)}_{N+n-2}) -a =0,
\end{equation*}
where $\ds a=e^{u}$, $\ds b= e^{\frac{n-2}{N+n-2}\xi}$ and $\ds \omega^{(n)}_{N+n-2} = e^{z^{(n)}_{N+n-2}\xi}$.
When $n=2$, we have the equation
$$ a(\omega^{(2)})^{3} - (1+a^{2})(\omega^{(2)})^{2} + (a^{2}+1)(\omega^{(2)}) -a =0 $$
By subtracting two equations we get
\begin{align*}
&a((\omega^{(n)}_{N})^{3}-(\omega^{(2)})^{3}) -(a^2 +1)((\omega^{(n)}_{N})^2 - (\omega^{(2)})^{2})+ (a^2 +1)(\omega^{(n)}_{N} -\omega^{(2)})  \\
=& -(ab^2 - a)(\omega^{(n)}_{N})^3 + (b^2 + a^2 b -a^2 - 1)(\omega^{(n)}_{N})^2 - (b-1) (\omega^{(n)}_{N})
\end{align*}
This implies
\begin{align*}
& \omega^{(n)}_{N} -\omega^{(2)} \\
=& \frac{-(ab^2 - a)(\omega^{(n)}_{N})^3 + (b^2 + a^2 b -a^2 - 1)(\omega^{(n)}_{N})^2 - (b-1) (\omega^{(n)}_{N})}{a((\omega^{(n)}_{N})^2 + \omega^{(n)}_{N}\omega^{(2)} + (\omega^{(2)})^2) - (a^2 +1)(\omega^{(n)}_{N} + \omega^{(2)}) +(a^2 +1) } \\
=& (b-1) \frac{-a(b+1)(\omega^{(n)}_{N})^3 + ( (b+1)+a^2)(\omega^{(n)}_{N})^2 - (\omega^{(n)}_{N})}{a((\omega^{(n)}_{N})^2 + \omega^{(n)}_{N}\omega^{(2)} + (\omega^{(2)})^2) - (a^2 +1)(\omega^{(n)}_{N} + \omega^{(2)}) +(a^2 +1) }
\end{align*}
For simplicity, we denote the right hand side by $(b-1)K^{(n)}_{N}$. Note that $ K^{(n)}_{N} \xrightarrow{N \to \infty} K \neq 0$.
On the other hand,
\begin{align*}
\omega^{(n)}_{N}  -\omega^{(2)}
= & e^{z^{(n)}_{N} \xi} - e^{z^{(2)} \xi} \\
= & e^{z^{(2)} \xi} ( e^{(z^{(n)}_{N} - z^{(2)}) \xi} - 1 ) \\
= & e^{z^{(2)} \xi} ( (z^{(n)}_{N} - z^{(2)}) \xi) ( \sum_{k=1}^{\infty} \frac{((z^{(n)}_{N} - z^{(2)}) \xi)^{k-1}}{k!})
\end{align*}

As a result,
\begin{align*}
&z^{(n)}_{N} - z^{(2)} \\
= &(b-1)\frac{K^{(n)}_{N}}{\xi e^{z^{(2)} \xi}( \sum_{k=1}^{\infty} \frac{((z^{(n)}_{N} - z^{(2)}) \xi)^{k-1}}{k!})} \\
= & \frac{n-2}{N+n-2} (\sum_{k=1}^{\infty} \frac{[(n-2)\xi /(N+n-2)]^{k-1} }{k!}) \frac{K^{(n)}_{N}}{ e^{z^{(2)} \xi}( \sum_{k=1}^{\infty} \frac{((z^{(n)}_{N} - z^{(2)}) \xi)^{k-1}}{k!})}\\
= & \frac{M^{(n)}_{N}}{N},
\end{align*}
where $M^{(n)}_{N} \xrightarrow{N \to \infty} M < \infty$. Therefore, we have
\begin{align*}
\lim_{N \to \infty} N(\Phi^{(2)}(z^{(n)}_{N}) - \Phi^{(2)}(z^{(2)}))
=\lim_{N \to \infty} \frac{\Phi^{(2)}(z^{(2)} + \frac{M^{(n)}_{N}}{N}) - \Phi^{(2)}(z^{(2)})}{\frac{M^{(n)}_{N}}{N}} M^{(n)}_{N}
=  0,
\end{align*}
where in the last equality we use the fact that $z^{(2)}$ is the solution of the saddle point equation 
$$\dfrac{d\Phi(z^{(2)})}{dz} = 0.$$
\end{proof}

\section{Evaluation at other root of unity}\label{sec4}

In this section we consider the behavior of $SU(n)$ invariant at other root of unity. Recall that conjecture~\ref{cvcn} is true for figure eight knot \cite{CLZ16}. Therefore it is natural to see whether our main theorem can be extended to other root of unity using the same tricks as before. Unfortunately the same trick does not apply at $\ds q=e^{\frac{2\pi i + u}{N+a}}$~when $a < n-2$.\\

To see why is it so, first, applying Lemma~\ref{Slemma} with the values
$$ \gamma=\frac{2\pi - iu}{2(N+a)}, \quad \xi = 2\pi i + u \quad \text{and}\quad z=\pi - iu -2(a+l)\gamma$$\\
and observing that $\ds \frac{\xi}{N+a}=2i\gamma $, we have
\begin{align}
\prod_{l=1}^{k}(1-e^{\frac{N-l}{N+a}\xi}) = \frac{S_{\gamma}(\pi - iu - (2(a+k)+1)\gamma)}{S_{\gamma}(\pi - iu - (2a+1)\gamma)} \label{G1}
\end{align}
Similarly, put $z=-\pi -iu +2(l+n-2-a)\gamma$, we have
\begin{align}
\prod_{l=1}^{k}(1-e^{\frac{N+l+n-2}{N+a}\xi}) = \frac{S_{\gamma}(-\pi - iu + (2(n-1-a) - 1)\gamma)}{S_{\gamma}(-\pi - iu + (2(n-2-a+k)+1)\gamma)}
\label{G2}
\end{align}
On the other hand,
\begin{align*}
\frac{[n-2+k]!}{[k]!}
&= [n-2+k][n-2+k-1]\dots [k+2][k+1]\\
&=q^{-\frac{n-2}{2}(\frac{n+2k-1}{2})}\prod_{l=k+1}^{n-2+k}(1-q^{l})
\end{align*}
Put $\ds q=\exp (\frac{\xi}{N+a})$, we have
\begin{align*}
1-q^{l}=1-e^{\frac{\xi l}{N+a}}
= 1 + e^{i(-\pi+2l\gamma)}
\end{align*}

By Lemma~\ref{Slemma} with value $z=-\pi  +2l\gamma$, we have
\begin{align}
\prod_{l=k+1}^{n-2+k}(1 + e^{i(-\pi+2l\gamma)})
&=\prod_{l=k+1}^{n-2+k} \frac{S_{\gamma}(-\pi  +(2l-1)\gamma)}{S_{\gamma}(-\pi + (2l+1)\gamma)}\notag \\
&=\frac{S_{\gamma}(-\pi + (2(k+1)-1)\gamma)}{S_{\gamma}(-\pi + (2(n+k-2)+1)\gamma)} \label{G3}
\end{align}

By (\ref{G1}), (\ref{G2}) and (\ref{G3}) we have
\begin{align*}
J_{N}^{(n)}(4_{1}, e^{\frac{\xi}{N+a}})
= &\frac{1}{[n-2]!}\frac{S_{\gamma}(-\pi-iu+(2(n-1-a)-1)\gamma)}{S_{\gamma}(\pi-iu-(2a+1)\gamma)} \times \\
&\sum_{k=0}^{N-1}[e^{-ku-(k(n-2-a)+(\frac{n-2}{2})(\frac{n-1}{2}))(\frac{\xi}{N+a})} \times\\
&  \frac{S_{\gamma}(\pi-iu-(2(a+k)+1)\gamma)}{S_{\gamma}(-\pi-iu+(2(n-2-a+k)+1)\gamma)}\frac{S_{\gamma}(-\pi+(2k+1)\gamma)}{S_{\gamma}(-\pi+(2n+2k-3)\gamma)}]
\end{align*}

Define
\begin{align*}
g_{N+a}(\omega) = &\exp(-(N+a)(u+n-2-a)\omega)\\
&\frac{S_{\gamma}(\pi - iu + i(\omega + \frac{a}{N+a})\xi)S_{\gamma}(-\pi - i\omega\xi)}{S_{\gamma}(-\pi - iu - i(\omega+\frac{n-2-a}{N+a})\xi)S_{\gamma}(-\pi - i(\omega+\frac{n-2}{N+a})\xi)}
\end{align*}

Since $S_{\gamma}(z)$ is defined for $|\operatorname{Re}(z)|<\pi + \operatorname{Re}(\gamma)$, one may check that $g(z)$ is well-defined when $z=x+iy \in D_{N+a}$, where $D_{N+a}$ is defined to be the set
$$D_{N+a}=
\left\{
 x+iy  | 
\begin{cases}
\quad- \frac{2\pi}{u} (x+\frac{a}{N+a}) -\frac{\operatorname{Re}(\gamma)}{u} < y < \frac{2\pi}{u} - \frac{2\pi}{u}(x+\frac{a}{N+a}) + \frac{\operatorname{Re}(\gamma)}{u} \\
\quad \quad \quad \quad \quad - \frac{2\pi x}{u} - \frac{\operatorname{Re}(\gamma)}{u}  < y < \frac{2\pi}{u} - \frac{2\pi}{u}x +  \frac{\operatorname{Re}(\gamma)}{u} \\
 - \frac{2\pi}{u} (x+\frac{n-2-a}{N+a}) -\frac{\operatorname{Re}(\gamma)}{u} < y < \frac{2\pi}{u} - \frac{2\pi}{u}(x+\frac{n-2-a}{N+a}) + \frac{\operatorname{Re}(\gamma)}{u} \\
\quad - \frac{2\pi}{u} (x+\frac{n-2}{N+a}) -\frac{\operatorname{Re}(\gamma)}{u} < y < \frac{2\pi}{u} - \frac{2\pi}{u}(x+\frac{n-2}{N+a}) + \frac{\operatorname{Re}(\gamma)}{u}
\end{cases}\right\}$$

We want to apply Residue theorem as in section~\ref{sec2}. To do so, we need to make sure that the poles of the function $z \mapsto \tan(N+a)\pi z$, namely $ \ds (\frac{2k+1}{2(N+a)},0)$ for $k=0,1,2, \dots, N-1$, sit inside the domain. It turns out that the poles sit inside the domain only when $a >n-3$. As a result, the proof cannot apply to the case where $\ds q=e^{\frac{2\pi i + u}{N+a}}$.

\section{Conclusion}\label{sec5}

In this section we discuss the difficulties about this project and give some remark about the proof.

First of all although there are already several results about the higher dimensional Reidemeister torsion, explicit values of the torsion is not known for most cases, in particular for figure eight knot. So we cannot compare our theorem with the exact value of higher dimensional Reidemeister torsion of the figure eight knot and draw any conclusion yet.

Nonetheless, if the conjecture is true, we can (i) obtain the higher dimensional Reidemeister torsion for hyperbolic knot explicitly by considering the asymptotic expansion of $SU(n)$ invariant and (ii) try to obtain some kind of relation (e.g. recursion formula) between Reidemeister torsion of different dimensions.

Another property revealed from the calculation is that the function $\Phi^{(n)}_{N}$ goes to $\Phi^{(2)}$ as $N$ goes to infinity. It is interesting to see whether it is true for other cases. In particular in \cite{HM14} H.Murakami consider the case where $K$ is twice-iterated torus knot. The authors hope to do this calculation in the future.

\section*{Acknowledgments}
We thank H.~Murakami for helpful comments in the preparation of this manuscript.  The work is done in the MPhil study of the first author, who received guidance from teachers and fellow students of the Department of Mathematics, The Chinese University of Hong Kong, in addition to various kinds of support including financially.  In particular, special thanks would be given to Prof.~K.W. Chan, Prof. Z.T. Wu, Marco Suen and Savio Chung.

\end{document}